\documentclass[ijoo,nonblindrev]{informs-ijoo}

\OneAndAHalfSpacedXI

\usepackage[square]{natbib}
 \bibpunct[, ]{(}{)}{,}{a}{}{,}%
 \def\bibfont{\small}%
\TheoremsNumberedThrough     
\ECRepeatTheorems

\EquationsNumberedThrough    

\MANUSCRIPTNO{}

\usepackage{appendix}

 \usepackage[font=small]{subcaption}
\usepackage{blkarray, bigstrut}
 \usepackage{graphicx}

 \usepackage{nomencl} 

\usepackage{pgfplots}
\usepgfplotslibrary{fillbetween}
\pgfplotsset{
  compat=newest,
  xlabel near ticks,
  ylabel near ticks
}

\usepackage{soul} 

 \usepackage[colorlinks=true, allcolors=blue]{hyperref}
 \usepackage{tikz}  
 \usetikzlibrary{calc} 
 \usepackage{circuitikz} 
 \usepackage{bm}
\usepackage{nomencl} 
\makenomenclature

\nomenclature[A,01]{$\mathcal{N}$}{Set of vertices where $\mathcal{N} = \{1, \dots, n\}$}
\nomenclature[A,02]{$\mathcal{C}$}{Set of components where $\mathcal{C} = \{1, \dots, n\}$}
\nomenclature[A,03]{$\mathcal{E}$}{Set of edges}
\nomenclature[P,04]{$\mathbf{A}$}{Adjacency matrix with elements 1 for the existing edges in a network and 0 otherwise}
\nomenclature[P,05]{$B^A$}{Budget for disabling vertices in the attack problem}
\nomenclature[P,05]{$B^R$}{Budget for link addition in the response problem}
\nomenclature[P,06]{$a_{i}$}{Cost of disabling vertex $i$ in the attack problem}
\nomenclature[P,06]{$d_{ij}$}{Cost of adding the edge between vertices $i$ and $j$ in the response problem}
\nomenclature[V,07]{$v^A_{ic}$/$v^R_{ic}$}{ Binary variable indicating that if vertex \textit{i} belongs to component \textit{c} in the attack/response problem}
\nomenclature[V,08]{$b^A_{c}$/$b^R_{c}$}{Binary variable indicating if component \textit{c} is non-empty in the attack/response problem}
\nomenclature[V,09]{$y^A_{ij}$}{Binary variable indicating 1 if vertex \textit{i} is connected to vertex \textit{j} in the attack problem}
\nomenclature[V,09]{$y^R_{ij}$}{Binary variable indicating 1 if an intra-component link between vertex \textit{i} and vertex \textit{j} is added in the response problem}
\nomenclature[V,10]{$q^R_{ij}$}{Binary variable indicating 1 if vertex \textit{i} and vertex \textit{j} are connected in the response problem}
\nomenclature[V,11]{$\alpha^A$/$\alpha^R$}{Length of the largest component in the attack/response problem}

\usepackage{todonotes}

\usepackage{interval}

\begin{document}



\RUNAUTHOR{Dehghani Filabadi and Bagheri}
\RUNTITLE{ Robust-and-Cheap Framework for Network Resilience} 
\TITLE{Robust-and-Cheap Framework for Network Resilience: A Novel Mixed-Integer Formulation and Solution Method}



\ARTICLEAUTHORS{%
\AUTHOR{Milad Dehghani Filabadi}
\AFF{Department of Integrated Systems Engineering, The Ohio State University, \EMAIL{dehghanifilabadi.1@osu.edu}} 
} 

\ABSTRACT{%
Resilience and robustness are important properties in the reliability and attack-tolerance analysis of networks. In recent decades, various qualitative and heuristic-based quantitative approaches have made significant contributions in addressing network resilience and robustness. However, the lack of exact methods such as mixed-integer programming (MIP) models is sensible in the literature. In this paper, we contribute to the literature on the network resilience and robustness for targeted and random attacks and propose a MIP model considering graph-theoretical aspects of networks. The proposed MIP model consists of two stages where in the first stage the worst-case attack 
is identified. Then, the second stage maximizes the network resilience under the worst-case attack by adding links considering a link addition financial budget. In addition, we propose a solution method 
that (i) provides a tight relaxation for the MIL formulation by relaxing some of the integrality restrictions, (ii) exploits the structure of the problem and reduces the second-stage problem to a less complex but equivalent problem, and (iii) identifies underlying knapsack constraints 
and generates lifted cover inequalities (LCI) cuts for identified knapsack constraints. We conclude numerical experiments for random networks and then extend our results to power system networks. Numerical experiments demonstrate the applicability and computational efficiency of the proposed \textit{robust-and-cheap} framework for network resilience. 
}


\maketitle

\section{Introduction} \label{introduction}
Resilience and robustness are important properties in the vulnerability and attack-tolerance analysis of networks such as the Internet, communication networks, transportation networks, flight networks, power systems, etc.
In case of a network failure due to targeted attacks, natural disasters, political unrest, human errors, equipment failures, etc., a network may divide into several components by losing some of its critical nodes and connections. The terms ``resilience” and ``robustness'' are often interchangeably used in the literature as improving one can lead to improving the other \citep{bachmann2020survey}. In particular, resilience is the capability of a system to maintain its functions and the ability to return to normal conditions after the occurrence of a disruptive event \citep{ahmadian2020quantitative}. Robustness measures the hardiness of a system against serious disruptions and refers to the ability of a system to remain functioning under a range of disturbances \citep{mens2011meaning}. 
There are several examples of the consequences of cascading failures in non-resilient networks. For example, the Japanese earthquake and tsunami in 2011 led to more than 15,000 deaths and disturbances within global supply chain networks \citep{mackenzie2012measuring}. In October 2012, Hurricane Sandy hit the U.S East Coast led to power outages to over 8 million customers across 21 states for a few weeks \citep{henry2016impacts}. More recently, the COVID-19 pandemic started in 2019 led to global disruptions in the supply chain network due to the lack of resilience in global supply chains \citep{golan2020trends}. In particular, at least 94\% of Fortune 1000 companies have encountered supply chain disruptions as a
result of the COVID-19 pandemic \citep{ivanov2020predicting}.
Such inevitable failures have captured increasing attention among researchers in recent decades. A promising solution in the literature is to enhance the network resilience and robustness by modifying the topological structure of networks to mitigate the impacts of such consequences in the face of failures or attacks to the parts of the network \citep{kazawa2020effectiveness}. 

Enhancing network resilience and robustness leads to questions about (i) resilience and robustness measures, (ii) topological structure of networks, and (ii) protection/modification methodologies to improve the network resilience and robustness \citep{bachmann2020survey}. In terms of definitions, the resilience and robustness measures proposed in the literature vary in different fields, e.g., mathematics, physics, computer science, and biology, to capture the connectedness of a network in some way. In particular, \cite{albert2000error} considered
the relative size of the largest connected component and the average size of the other components, \cite{holme2002attack} and \cite{beygelzimer2005improving} considered mean shortest paths, \cite{latora2007measure} and \cite{scellato2011evaluating} used efficiency-based measures to quantify the shortest path distances between pairs of stations (nodes) in the network. The fraction of connected node pairs was studied by \cite{sun2007pairwise}, and \cite{matisziw2009modeling} defined a metric based on the available flows between the node pairs. A comprehensive survey on robustness analysis and robustness measures is provided in \cite{bachmann2020survey}.

From the topological perspective of networks, several studies addressed network resilience and robustness based on graph-theoretical properties. In particular, such studies measured resilience and robustness by considering three important quantities, namely: (1) the number of nodes that are not functioning/attacked, (2) the number of remaining connected subnetworks, and (3) the size of the largest remaining component within which the nodes are still connected. Examples of such measures include toughness \citep{chvatal1973tough}, scattering degree \citep{jung1978class}, integrity \citep{barefoot1987vulnerability}, tenacity \citep{cozzens1995tenacity}, and rupture degree \citep{li2005rupture}. Among such graph-theoretical measures, only the last two parameters, i.e., tenacity and rupture degree, consider all three important quantities of the network. The rupture degree is the newest and is a more comprehensive measure that represents a trade-off between the amount of required work to damage the network and the severeness of the damage in the network \cite{li2010extremal}.

From the methodological point of view, several protection/modification strategies have been proposed in the literature to improve network resilience and robustness. The existing strategies can be categorized into three main groups:  The first group is to protect the critical nodes. For example, \cite{ruj2014analyzing} showed that cascading failures are mainly due to targeted attacks, rather than random attacks, and proposed a strategy to preserves the critical nodes of a network. \cite{wang2020multiple} proposes a node protection method based on multiple topology features. \cite{nguyen2013detecting} showed the critical node identification is an NP-hard problem and used a greedy framework to find such nodes. \cite{abbas2017improving} proposed a framework to identify a small subset of nodes that always remain intact, called trusted nodes. 
\cite{sen2014identification} proposed a heuristic framework to identify and protect the $K$ most critical nodes in an interdependent network. The second group of network resilience and robustness studies focus on adjusting the properties of the network. In particular,
\cite{zhao2015cascading} proposed to increase node capacities to improve the robustness, \cite{shao2011cascade} proposed a strategy to change the deployment of independent
nodes, and \cite{chan2016optimizing} 
modified a network by re-writing, or swapping, the existing edges while the degree of nodes remains unchanged. Such approaches are mainly practical in networks with no or limited physical infrastructures, such as the Internet, since it requires the existing properties of the network to be changed substantially.

The third group of studies focuses on adding new connecting links to improve the robustness of networks. There are several approaches for link addition such as random addition \citep{cao2013improving}, low-degree addition \citep{zhao2009enhancing},
low betweenness addition \cite{guan2011routing},
and algebraic connectivity-based addition \cite{wang2008algebraic}. In \cite{ji2016improving}, the authors studied two link addition methods for improving network robustness for random attacks. In \cite{cui2017enhance}, authors proposed to add links based on receiving capability ratio of nodes. \cite{dong2015approaches} proposed a method to add links based on the degree and betweenness. \cite{wang2018improving} proposed to add links based on neighbor node priority connection to improve robustness. In \cite{dong2020improving}, the authors proposed to add a larger number of short links dispersedly rather than adding a smaller number of long links concentratedly.

Link addition strategy is the more practical and effective method, among other alternatives, to improve the resilience and robustness of networks \citep{kazawa2020effectiveness}. 
However, the main limitation of the existing link addition strategies is that these studies are mainly based on heuristic approaches. Although such approaches improve the network resilience and robustness by finding a local optimal solution, there is no guarantee to find the best improvement in the resilience and robustness, i.e., global optimal solution. Therefore, the lack of exact optimization models, especially in particular mixed-integer programming (MIP) models, is sensible in the literature. For example, in network analysis, there are many binary decisions such as determining if a node is critical, if a link is impacted by attacking the critical nodes, or if adding a given link to the network improves its resilience and robustness, etc. Such decisions can be incorporated using MIP that allows variables to take integer values, such as 0 and 1 for binary decisions. 

MIP models can be solved using exhaustive searching through branch-and-bound (B\&B) algorithm which is computationally infeasible in large-scale problems due to the NP-completeness of MIP models \citep{kannan1978computational, wang2008algebraic}. Alternatively, branch-and-cut (B\&C) algorithms combine cutting plane algorithms with the B\&B algorithm to solve MIPs more time-efficiently \citep{wolsey1999integer}. Cutting plane algorithms, such as Gomory’s mixed-integer cuts \citep{gomory1960algorithm}, mixed-integer rounding inequalities \citep{nemhauser1990recursive}, simple and lifted knapsack cover cuts \citep{wolsey1975faces,balas1975facets}, intersection cuts \citep{balas1971intersection}, etc., reduce the computational challenges by generating valid linear inequalities that preserve all the integer solutions but potentially eliminate a part of the feasible region with fractional solutions. In particular, \cite{bixby2004mixed} empirically showed that when implementing cutting planes in commercial MIP solvers, the computational time significantly speeds up. In addition to cutting plane algorithms, reduction methods are efficient to decrease the complexity of a problem, including MIP problems. Such methods usually exploit the structure of a problem and reduce a given problem into a simpler problem with an equivalent optimal solution. Depending on the optimization model, several reduction methods have been proposed in the literature including \cite{guignard1981logical,cook1993implementation,tits2006constraint}. More related to our work, relaxation techniques \citep{geoffrion1974lagrangean,bsc2005relaxation} are among reduction methods where some of the restrictions on the problem, such as integrality constraints, are relaxed. The power of such relaxations depends on how "tight" they are where tighter relaxations correspond to relaxations whose optimal solutions are more close to the optimal solutions of the original problem.

In this paper, we contribute to the literature on network resilience and robustness against targeted or random attacks. We propose a framework to optimally increase the network hardiness against attacks, i.e., network robustness, and to increase the network functionality after attacks, i.e., network resilience. To the best of our knowledge, we propose the first MIP formulation for network resilience and robustness considering graph-theoretical aspects of a network by modeling three important topological quantities of a network namely: (i) the number of attacked or failed elements, (ii) the number of remaining connected subnetworks after an attack occurs, and (iii) the size of the largest remaining group within which mutual communication can still occur after an attack. We first propose a two-stage MIP formulation for targeted attacks where in the first stage the worst-case attack leading to the most severe damage in the network is identified, and used for measuring the network robustness. In the second stage, after the worst-case attack is realized, we optimally improve the network resilience by adding links considering a financial budget constraint for link addition. Finally, we show the applicability of our proposed approach to other types of attacks. In addition, we propose a solution method based on cutting planes and relaxation methods to efficiently solve the proposed MIP model. Particularly, we provide a tight relaxation for the first-stage problem by relaxing some of the integrality constraints. Then, we exploit the structure of the problem and reduce the second-stage problem into a significantly less complex but equivalent problem. Finally, we identify underlying knapsack constraints in our formulation through a proposed re-indexing procedure and then manually implement knapsack lifted cover inequalities (LCI), which is known as a recent cutting plane algorithm. The main contributions of this paper are as follows:

\begin{enumerate}
\item We propose a novel MIP formulation for optimally improving network resilience and robustness. The proposed approach can be applied to targeted attacks through a two-stage problem and also can be used for random or designated attacks through a single-stage problem.
\item We propose a solution method to solve the proposed MIP formulation in a time-efficient manner. Particularly, the proposed solution method (i) provides a tight relaxation by removing some of the integrality restrictions, (ii) exploits the structure of the problem and reduces the second stage problem to a less complex but equivalent problem, and (iii) identifies underlying knapsack constraints through a re-indexing procedure and generates efficient cuts for knapsack constraints.
\item We implement the proposed solution method in several instances and demonstrate that our approach provides a \textit{robust-and-cheap} framework for network resilience problems due to the resulting computational efficiency and robustness. Moreover, we show how our proposed framework can provide a trade-off for the budget against network resilience and robustness.
\end{enumerate}

The rest of this paper is organized as follows: Section~\ref{section: definitions} presents the required definitions used in this paper. In Section~\ref{section : model}, the proposed two-stage optimization model for targeted attacks is presented and discussed in detail. Section~\ref{section: solution} proposes a solution method to efficiently solve the proposed MIP formulation. In Section \ref{section other attacks}, we show how the proposed framework is extended to other types of attacks. Section~\ref{section; results} provides numerical results and analyses for random networks and extends the results for power system application in Section \ref{section : power system extention}. Finally, concluding remarks are provided in Section~\ref{section: conclusion}.

Throughout this paper, paired terms (edge, link), (vertex, node), and (graph, network) might be interchangeably used.

\section{Definitions} \label{section: definitions}
In this section, we first present graph-theoretical definitions used in this paper. In the second part, we define the network resilience and robustness of networks and present metrics to measure such quantities.

\subsection{Graph Theory Background}
\begin{definition}
\textit{Graph} $G$ is the mathematical representation of a network and describes the relationship between vertex set $V(G)$ and edge set $E(G)$.
\end{definition}

\begin{definition}
A set $X \subseteq V(G)$ of graph $G$ is called a \textit{cut set} if $G-X$ is disconnected or $G-X$ has only one vertex, where term $G-X$ is the remaining graph obtained by removing vertices (and the corresponding edges) of $X$. 
\end{definition}
\begin{definition}\label{def connected graph}
A \textit{connected graph} is a graph where there exists at least one path between any two vertices of the graph. A graph is said to be \textit{disconnected} if it is not connected.
\end{definition}
\begin{definition} \label{def component} 
Each disconnected graph consists of two or more connected subgraphs that are called \textit{components} of the graph. The smallest non-empty component in a graph is an isolated vertex.
\end{definition}

\begin{definition}
Graph $G+e_{ij}$ is a graph obtained by adding edge $i$-$j$ to graph $G$.
\end{definition}
\begin{definition}\label{def rupture par}
\citep{li2005rupture} The \textit{rupture degree} of a connected graph $G$ is defined by
\begin{equation}\label{rupture}
r^*(G) = \max_{X} \{-|X|-m(G-X)+\omega(G-X)\}
\end{equation}
where the maximum is taken over all cut sets $X \subseteq V(G)$, and $\omega(G-X)$ and $m(G-X)$ denote the number of components and the number of vertices (length) of the largest component of graph $G-X$, respectively. 
\end{definition}
\begin{definition} \label{worst cut set}
The \textit{worst cut set} $X^*$ is the cut set determining to the rupture degree of a graph, i.e., 
\begin{equation}\label{eq worst-case cutset}
X^*=\text{arg}\max_{X} \{-|X|-m(G-X)+\omega(G-X)\}
\end{equation} 
\end{definition}


\subsection{Resilience and Robustness Definition}
One important factor in resilient networks is to respond to an attack so that the network maintains its functioning after the attack. There are two factors associated with each attack: (i) how hard it is to damage a network, i.e., the network robustness, and (ii) how hard it is to maintain the proper functioning of the remaining parts of the network after an attack occurs, i.e., the network resilience. In the context of network resilience and robustness, we use the definition of the rupture degree and interpret the worst cut set as the attack leading to the most severe damage in the network, i.e., the \textit{worst-case attack}. We quantify the network robustness by the size of the worst cut set ($|X^*|$ from Definition \ref{worst cut set}) since it corresponds to the network hardiness against the most severe attacks. Also, we measure the network resilience by the negative of the rupture degree ($-r^*(G)$ from Definition \ref{def rupture par}) since it determines how badly a network is damaged by the worst-case attack, or equivalently how functioning the network is after the worst-case attack. For example, a network with smaller cut sets can be damaged easier (less robust) because a fewer number of vertices should be disabled to most severely damage the network. Also, a graph with a higher rupture degree is less functioning after the attack (less resilient) because it is more badly damaged and ruptured by the attack. However, as the rupture degree decreases, the network becomes more resilient since the underlying graph has stronger connectivity properties.

\section{Optimization Model}\label{section : model}
To address both the resilience and robustness concerns of a network, one can find the worst-case attack by finding the most critical vertices (i.e., the worst cut set) of a network. Then, in the absence of the most critical vertices that are supposedly disabled by the attack, proper actions such as link addition can be carried out to maximize the network, i.e., minimizing the rupture degree, so that the network is resilient (immunized) against the worst-case attack. 

Optimizing under the worst-case scenario is a widely-used approach to immunize the system not only under the worst-case scenario but also under any other possible scenarios \citep{bertsimas2004price, filabadi2020robust, filabadi2019effective}. In the context of network resilience and robustness, optimizing the network resilience and robustness is obtained by adding new links under the worst-case attack so that the \textit{reconstructed network}, i.e., the initial network with new links, is less dependant on the most critical vertices disabled by the worst-case attack. As a result, since the network resilience and robustness are improved under the worst-case attack, the resilient network is more robust and functional under any other attacks which are likely to occur.

To model the network resilience and robustness, we propose a two-stage optimization model of form 
\begin{subequations}\label{two stage general}
    \begin{align} 
\min_{\mathbb{E}^R} \quad \Big( f^R({\mathbb{E}^R,\mathbb{E}^{A*}}) \quad +\quad & \max_{\mathbb{E}^A}  \quad  f^A(\mathbb{E}^A) \Big) \label{obj func} \\ 
& \quad \text{s.t.} \quad h^A(\mathbb{E}^A) \le l^A \label{stage 1 const}\\
\text{s.t.} \quad h^R(\mathbb{E}^R,\mathbb{E}^{A*}) \le l^R \label{stage 2 const}
\end{align}
\end{subequations}
where the first stage, called \textit{attack problem} and denoted by superscript $A$, consists of the inner minimization problem along with constraint \eqref{stage 1 const} to find the worst-case attack, i.e., the attack corresponding to the maximum value of the rupture degree $f^A$ (or equivalently minimum value of the network resilience). Given the worst-case attack (i.e., $\mathbb{E}^{A*}=\text{arg}\max_{\mathbb{E}^A} f^A(\mathbb{E}^A)$) is realized, the outer minimization problem, called \textit{response problem} and denoted by superscript $R$, aims to minimize the rupture degree $f^R$ (or equivalently maximize the network resilience) by adding links under the worst-case attack. Constraint set \eqref{stage 2 const} corresponds to structural constraints and financial budget constraints for link addition. The proposed problem is inherently a bilevel problem where the lower level problem (attack) problem is embedded in the upper level (response) problem. However, since the decision variable sets $\mathbb{E}^A$ and $\mathbb{E}^R$ are distinct and are determined independently in a sequential manner, formulation \ref{two stage general} can be decomposed in two separate problems where we first solve the attack problem \ref{stage 1} and, then given optimal solution $\mathbb{E}^{A*}$, we solve the response problem. The output of the response problem is a \textit{reconstructed} network consisting of the initial network and additional links. 

In what follows, we first present each stage separately. Then in later sections, we propose a solution method to solve the optimization problem. The following notations are used in the proposed formulation of this paper, where we use boldface uppercase and lowercase characters to denote matrices and vectors, respectively.

\printnomenclature

\subsection{Stage I: Attack Problem}
This stage aims to find the most critical vertices by finding the worst-case attack that can lead to the most severe damages in a network (maximum value of the rupture degree). Once a vertex is attacked, the vertex and all its connections are removed from the network. Given the set of decision variables $\mathbb{E}^A := \{ \mathbf{v}^A, \mathbf{b}^A, \mathbf{y}^A, \alpha^A \}$, the objective function of the attack problem is defined in \eqref{s1 obj func}, based on the rupture parameter, and consists of three parts. The first part finds the number of disabled vertices which is the difference of all vertices of the network, i.e., $|\mathcal{N}|$, and the active vertices under the worst-case disruption, i.e., $\sum_{i \in \mathcal{N}} \sum_{c \in \mathcal{C}} v^A_{ic}$. The last two terms in \eqref{s1 obj func} correspond to the length of the largest component and the number of components, respectively.
\begin{equation}
\max_{\mathbb{E}^A}  \quad  -(|\mathcal{N}|-\sum_{i \in \mathcal{N}} \sum_{c \in \mathcal{C}} v^A_{ic}) - \alpha^A + \sum_{c \in \mathcal{C}} b^A_c  \label{s1 obj func} \tag{4a}
\end{equation}

Constraints \eqref{s1 c1}-\eqref{s1 c8} address the operational, physical, and financial limits in the network such as connectivity and budget constraints.
\begin{subequations}\label{stage 1}
\begin{align} 
& \sum_{c \in \mathcal{C}} v^A_{ic} \le 1 & \forall i \in \mathcal{N} \label{s1 c1} \tag{4b} \\
& \sum_{i \in \mathcal{N}} v^A_{ic} \le \alpha^A & \forall c \in \mathcal{C} \label{s1 c2} \tag{4c} \\
& b^A_c \le \sum_{i \in \mathcal{N}} v^A_{ic} & \forall c \in \mathcal{C} \label{s1 c3} \tag{4d} \\
& \sum_{\substack{ j \in \mathcal{N} \\ i \neq j}} y^A_{ij} \le (|\mathcal{N}|-1)\sum_{c \in \mathcal{C}} v^A_{ic} & \forall i \in \mathcal{N} \label{s1 c4} \tag{4e} \\
& \sum_{i \in \mathcal{N}} a_i \Big(1-\sum_{c \in \mathcal{C}} v^A_{ic} \Big) \le B^A  \label{s1 c5} \tag{4f} \\
& y^A_{ij} = y^A_{ji} & \forall i,j \in \mathcal{N} \label{s1 c6} \tag{4g} \\ 
&A_{ij}(\sum_{c \in \mathcal{C}} v^A_{ic}+\sum_{c \in \mathcal{C}}v^A_{jc}-1) \le y^A_{ij} & \forall i,j \in \mathcal{N}: i>j \label{s1 c6.2} \tag{4h} \\
& y^A_{ij} \le A_{ij}(1-v^A_{ic}+v^A_{jc}) & \forall i,j \in \mathcal{N}: i>j, \forall c \in \mathcal{C} \label{s1 c8} \tag{4i} \\
& y^A_{ij} \le A_{ij}(1+v^A_{ic}-v^A_{jc}) & \forall i,j \in \mathcal{N}: i>j, \forall c \in \mathcal{C} \label{s1 c9} \tag{4j} \\
&v^A_{ic} \in \{0,1\}, y^A_{ij} \in \{0,1\}, b^A_{c} \in \{0,1\}, \alpha \in \mathbb{Z}_+ \in \{0,1\}& \forall i,j \in \mathcal{N}, \forall c \in \mathcal{C} \label{domain variable s1} \tag{4k}
\end{align}
\end{subequations}



Constraint \eqref{s1 c1} ensures that each vertex is at most in one component, if active, or it is disabled. Given $\sum_{i \in \mathcal{N}} v^A_{ic}$ is the length of component $c$, constraint \eqref{s1 c2} enforces the length of the largest component $\alpha$ be equal or greater than the length of all components. Constraint \eqref{s1 c3} ensures that if the length of a component is zero, the component is empty. Otherwise, due to the direction of the objective function, $b_c^A=1$ indicating the component is non-empty. Constraint \eqref{s1 c4} shows that each vertex $i$ is connected to at most $|\mathcal{N}|-1$ vertices, if the vertex is active under the worst-case attack. Otherwise, all edges connected to vertex $i$ are eliminated. It is noteworthy that $y_{ii}$ is excluded from the left-hand size expression in \eqref{s1 c4} since there does not exist an edge with length 0 in the network. Constraint \eqref{s1 c5} is the budget constraint for attacking the most critical vertices. Constraint \eqref{s1 c6} shows that the network is undirected. Constraint \eqref{s1 c6.2} ensures that if two vertices $i$ and $j$ are connected in the initial network, i.e. $A_{ij}=1$, and they remain in the same component $c$ after attack (i.e., $v^A_{ic}=v^A_{jc}=1$), these vertices remain connected. In other words, this constraint preserves the structure of parts of the network whose connections are not impacted by the worst-case attack.
Constraints \eqref{s1 c8}-\eqref{s1 c9} are complementary constraints that correspond to the connection of two vertices depending on their position. For instance, 
if two active vertices $i$ and $j$ are in different components, Constraints \eqref{s1 c8} and \eqref{s1 c9} reduce to $y^A_{ij}\le \min\{0,2A_{ij}\}$, enforcing that the two vertices are not connected. 
The variable domains are represented in constraint \eqref{domain variable s1}.

It is to be mentioned that the number of components is not known in advance since it depends on the structure of the network after removing the most critical vertices. However, the maximum number of components is $|\mathcal{C}|=|\mathcal{N}|$ which is the number of vertices in the graph since an isolated vertex corresponds to the smallest non-empty component. If for a component $c$ we have $\sum_{i \in \mathcal{N}} v^A_{ic}=0$, it means component $c$ is empty and it is not a part of the remaining network after the attack.

If there is enough budget for attacking, the outcome of the attack problem is a disconnected network with the highest value for the rupture degree. This solution can be used in the response problem to optimize network resilience after the attack.

\subsection{Stage II: Response Problem}
The link addition strategy is an effective method to improve the robustness of networks \citep{kazawa2020effectiveness}. The response problem in the second stage is carried out after the realization of the worst-case attack where its inputs are the optimal solution of the attack problem, i.e., $\mathbb{E}^{A*} := \{ \mathbf{v}^{A^*}, \mathbf{b}^{A^*}, \mathbf{y}^{A^*}, \alpha^{A^*} \}$, to proactively construct a resilient network that is immunized against the worst possible attacks. Particularly, after the critical nodes corresponding to the worst-case attack are identified and removed, the response problem aims to maximize the network resilience under the worst-case attack while considering a budget for link addition. Therefore, adding links under the worst-case attack ensures that the reconstructed network is able to maintain its function not only in the absence of the most critical nodes but also in the absence of other sets of nodes.

In the response problem, the available vertices are those that remain active (not attacked) after the worst-case attack. The objective function aims to find a reconstructed network corresponding to the minimum rupture degree. This can be obtained by adding links to the initial network that can possibly reduce the number of components and/or increase the length of the largest component $\alpha$ (definition \ref{def rupture par}). In what follows, we first show adding some candidate links is ineffective and would not change the rupture degree, or equivalently. Then, we identify only those links that can increase network resilience and provide a MIP formulation for link addition.

Given the attacked network is a disconnected network, there are two sets of candidates for link addition: (i) links within each component, called \textit{local links}, or (ii) links between components, called \textit{inter-component links}, to connect isolated components.
Adding local links within a component does not impact the rupture parameter, and consequently the objective function of the response problem. It is because each component is a connected subgraph with $V$ vertices, and the local links do not connect the current component to any new vertices, implying that the length of the largest component remains unchanged. Also, adding local links does not impact the number of components since adding such links just changes the connectivity structure of each individual component. Therefore, after removing the worst cut set $X$ from graph $G$, we need to add inter-component links to $G-X$ to decrease the rupture parameter and improve the network resilience.
Given the definitions, We first propose a lemma to establish a connection between graph $G$ and the length of its largest component, denoted by $m(G)$. Then we provide a proposition showing the inter-component links will result in a lower value for the rupture parameter.

\begin{lemma} \label{lemma 1}
Let $G$ be a graph and $\mathcal{C}$ be the set of its component. Then $m(G)=m(c_i)$ if and only if $c_i$ is the largest component of graph $G$.
\end{lemma}
\proof{Proof. } Given $m(G)=m(c_i)$, from the definition of a component it follows that the maximum number of connected vertices in graph $G$ equals to the number of vertices in component $c_i$, implying $c_i$ is the largest component. On the other hand, if $c_i$ is the largest component of $G$, the second direction simply follows from definitions. \Halmos

\begin{proposition} \label{prop 1}
Let $G-X$ be a disconnected graph with components $\mathcal{C}= \{c_1, \dots , c_s\}$. Given a link $e_{ij}$ connecting vertices $i$ and $j$, 
\begin{enumerate}
    \item[a.] \begin{equation} \label{prop 1 part a}
r(G-X+e_{ij}) \le r(G-X) - 1
\end{equation}
if $i \in c_l$, $j \in c_k$, and $l \neq k$.

\item[b.] In addition to (a), let $L = m(c_l)$, and $K=m(c_k)$. If $L+K > m(G-X)$ then 
\begin{equation} \label{prop 1 part b}
r(G-X+e_{ij}) \le r(G-X)-1 - \big(L+K-m(G-X)\big)
\end{equation}
\end{enumerate}
\end{proposition}

\proof{Proof.}
$a$. Let $G-X$ be a disconnected graph after removal of cut set $X$. Also let $c_l$ and $c_k$ be two distinct and arbitrary components of $G-X$ where $m(c_l)=L$ and $m(c_k)=K$.
Since each component $c_l$ and $c_k$ are connected subgraphs (from definition \ref{def component}), adding edge $e_{ij}$ unifies two components $c_l$ and $c_k$ and results in a larger component $c_l + c_k + e_{ij}$, and reduces the number of components by 1. Thus, $\omega(G -X + e_{ij}) = \omega(G-X)-1$, which implies $r(G-X+e_{ij}) \le r(G-X)-1$.\\
\noindent $b$. In addition to the conditions of a, assume $L+K > m(G-X)$. Therefore, adding $e_{ij}$ to graph $G-X$ results in a larger component $c_l + c_k + e_{ij}$ so that $m(c_l + c_k + e_{ij})=L+K>m(G-X)$. Thus, component $c_l + c_k + e_{ij}$ becomes the largest component of graph $G-X+ e_{ij}$ and from Lemma \ref{lemma 1}, $m(G-X+ e_{ij})=m(c_l + c_k + e_{ij})=L+K$, implying the length of the largest component has increases by $L+K-m(G-X)$. Thus, from this point and also part (a) of the proposition, $r(G-X+e_{ij}) \le r(G-X)-1 -(L+K-m(G-X))$. \Halmos

Proposition \ref{prop 1} shows how adding inter-component links guarantees improving the network resilience. Since inter-component links are the only set of links improving the network resilience, we allow only the addition of intra-component links in the response problem. 

After the optimal solution of the attack problem is obtained, one can easily identify the set of remaining components ($\mathcal{C}^R$), the set of remaining vertices ($\mathcal{C}^R$), and the set of vertices in each components ($V(c), c \in \mathcal{C}^R$). Let $t^R_c \in \{0,1\}, \forall c \in \mathcal{C}^R$ indicating 1 if component $c$ is the largest component after adding links in the response problem. Thus, given decision variables $\mathbb{E}^R := \{ \mathbf{v}^R, \mathbf{b}^R, \mathbf{y}^R, \boldsymbol{\alpha}^R, \mathbf{t}^R \}$, the objective function \eqref{s2 obj func} aims to maximize the network resilience by minimizing the rupture degree.
\begin{subequations} \label{response problem}
\begin{align}
-|X| \hspace{0.2cm}+\hspace{0.2cm} \min_{\mathbb{E}^R}  \quad & \Big( -\alpha^R + \sum_{c \in \mathcal{C}^R} b^R_c + \epsilon \sum_{c \in \mathcal{C}^R} t^R_c \Big) \label{s2 obj func} \\
&\sum_{i \in \mathcal{N}^R} v^R_{ic} \le \alpha^R \le \sum_{i \in \mathcal{N}^R} v^R_{ic} + M(1-t^R_c)  \qquad \qquad \forall c \in \mathcal{C}^R    \label{s2 aux1} \\
& \sum_{c \in \mathcal{C}^R} t^R_c \ge 1 \label{s2 aux2}
\end{align}
\end{subequations}
where $|X| = (|\mathcal{N}|-\sum_{i \in \mathcal{N}} \sum_{c \in \mathcal{C}} v^{A*}_{ic})$ is the size of the worst cut set which was identified in the attack problem. Parameters $\epsilon$ and $M$ are pre-defined small and big numbers, respectively. Since the size of the worst cut set is a constant, and known, in the response problem, we have excluded it from the objective function, but it will be considered for evaluating the rupture degree. The first and second terms of the objective function correspond to the length of the largest component and the number of components after adding links, respectively. The third term of the objective function is a small penalty added to the objective function to ensure minimizing the rupture degree is dominating while addressing constraints corresponding to the largest component(s) (i.e., \eqref{s2 aux1} and \eqref{s2 aux2}). In particular, Constraint \eqref{s2 aux1} finds the length of the largest component. If component $c$ is the largest component then $t^R_c=1$ and the constraint is binding. Otherwise, due to the penalty considered in the objective function, $t^R_c=0$ and thus the right-hand side of constraint \eqref{s2 aux1} is a big number, meaning that component $c$ does not determine the length of the largest component, i.e., $\alpha^R$, as it is not the largest component. Constraint \eqref{s2 aux2} ensures that at least one component is the largest component.  

In the response problem, when links are added to the attacked network, the topology of the reconstructed network may change due to changing the set of vertices in each component, the length of the largest component, the number of components, etc. Constraints \eqref{s2 c1} to \eqref{domain variable} address the structural and budget constraints of the response problem.
\begin{subequations}\label{stage 2}
\begin{align} 
& \sum_{c \in \mathcal{C}^R} v^R_{ic} = 1 & \forall i \in \mathcal{N}^R \label{s2 c1} \tag{7d} \\
& \sum_{i \in \mathcal{N}^R} v^R_{ic} \le  |\mathcal{N}^R|b^R_c&  \forall c \in \mathcal{C}^R \label{s2 c3}  \tag{7e} \\
&\sum_{i \in \mathcal{N}^R} \sum_{\substack{j \in \mathcal{N}^R \\ j > i}} d_{ij}y^R_{ij} \le B^R \label{budget} \tag{7f}\\
& y_{ij}^R = y_{ji}^R & \forall i,j \in \mathcal{N}^R \label{s2 c6} \tag{7g} \\
& q_{ij}^R = q_{ji}^R & \forall i,j \in \mathcal{N}^R \label{s2 c6.1}  \tag{7h} \\
& q^R_{ij} \le \sum_{i \in V(c)}\sum_{j \in V(c')} y^R_{ij} \le Mq^R_{ij} & \forall i \in V(c), \forall j \in V(c'),  \forall c,c' \in \mathcal{C}^R: c\neq c' \label{eq1} \tag{7i}\\
& v^R_{ic}+v^R_{jc}-1 \le q^R_{ij} & \forall i \in V(c), \forall j \in V(c'),  \forall c,c' \in \mathcal{C}^R: c\neq c'  \label{eq2} \tag{7j} \\
& q^R_{ij} \le 1-v^R_{ic}+v^R_{jc} & \forall i \in V(c), \forall j \in V(c'),  \forall c,c' \in \mathcal{C}^R: c\neq c'  \label{eq3} \tag{7k} \\
& q^R_{ij} \le 1+v^R_{ic}-v^R_{jc} & \forall i \in V(c), \forall j \in V(c'),  \forall c,c' \in \mathcal{C}^R: c\neq c'  \label{eq4} \tag{7l} \\
& \mathbf{v}^R \in \{0,1\}^{|\mathcal{N}^R||\mathcal{C}^R|}, \mathbf{y}^R \in \{0,1\}^{|\mathcal{N}^R|^2}, \mathbf{q}^R \in \{0,1\}^{|\mathcal{N}^R|^2}, \mathbf{b}^R \in \{0,1\}^{|\mathcal{C}^R|}, \span \label{domain variable} \nonumber \\
& & \mathbf{t}^R \in \{0,1\}^{|\mathcal{C}^R|}, \alpha^R \in \mathbb{Z}_+ \tag{7m}
\end{align}
\end{subequations}
Constraints \eqref{s2 c1} guarantees that after we add links, each active vertex is exactly in one of the components of the reconstructed network.
Constraint \eqref{s2 c3} identifies non-empty components so that if the length of a component is greater than 0 after adding links, the component is non-empty (i.e., $b^R_c=1$). Otherwise, due to the direction of the objective function $b^R_c=0$ implying the component is empty. 
Constraint \eqref{budget} is the budget constraint limiting the total link added in the reconstructed network. 
Constraints \eqref{s2 c6} and \eqref{s2 c6.1} correspond to the properties of undirected networks.
Constraint \eqref{eq1} shows that if two components are connected, there is at least one intra-component link between the two components.
Constraint \eqref{eq2} enforces that if two vertices are in the same component, then there is a path between the two vertices, i.e., they are connected. Constraints \eqref{eq3} and \eqref{eq4} are complementary and ensure that if two vertices are in different components, they are not connected. Finally, constraint \eqref{domain variable} represents the variable domains for all variables.

Once the response problem is solved, the optimal intra-component links $\mathbf{y}^R$ are determined to be added to the initial network. Thus, the resulting reconstructed is the most resilient network against any attacks including the previously-identified worst-case attack given the budget. We will elaborate more on this point in the numerical results.









 

\section{Solution Methodology}\label{section: solution}
In the proposed formulation, both attack and response problems are integer programming (IP) models which are difficult to solve when the problem size increases, due to the NP-completeness of integer programming \citep{papadimitriou1981complexity}. Therefore, in what follows, we propose a procedure to solve each IP model more efficiently. Particularly, for the attack problem, we identify a knapsack constraint and implement a more recent procedure for generating knapsack cover inequalities (CI) \citep{letchford2019lifted}, which has not been implemented in the commercial solvers yet. We also provide a tight relaxation by relaxing some of the integrality constraints of the attack problem. Then, for solving the response problem, we derive an equivalent but smaller problem to the response problem and also identify an underlying knapsack problem for which we can provide CIs. To do so, we first provide preliminary definitions of IP and knapsack CIs that will be used throughout this section.
%
%
\subsection{Preliminary Definitions}
\begin{definition} \label{conv combination}
A point $x$ in $\mathbb{R}^n$ is \textit{convex combination} of points $x_1, \dots, x_q \in \mathbb{R}^n$, if $\exists \lambda_1, \dots, \lambda_q \ge 0$ such that
\begin{equation*}
    x = \sum_{j=1}^q \lambda_j x_j, \hspace{1cm}  \sum_{j=1}^q \lambda_j= 1 
\end{equation*}
\end{definition}
\begin{definition} \label{conv set}
A set $S \in \mathbb{R}^n$ is a \textit{convex set} if 
for any two points $x,y \in S$, the line segment connecting $x$ and $y$ is contained in $S$.
\end{definition}
\begin{definition} \label{conv hull}
Given a set $S \in \mathbb{R}^n$, the \textit{convex hull} of $S$, denoted by $\text{conv}(S)$, is the smallest convex set containing $S$. That is
\begin{equation*}
\text{conv}(S)= \{\sum_{j=1}^q \lambda_j x_j: \quad x_1, \dots, x_q \in S, \hspace{0.1cm} \lambda_1, \dots, \lambda_q \ge 0, \sum_{j=1}^q \lambda_j= 1 \}  \hspace{1cm}   
\end{equation*}
\end{definition}
\begin{definition}\label{face of polyhedron}
Let $P$ be a non-empty polyhedron in $\mathbb{R}^n$. A \textit{face} of a polyhedron $P$ is a set of form 
\begin{equation*}
    F:= P \hspace{0.1cm} \cap \hspace{0.1cm} \{x \in \mathbb{R}^n \hspace{0.1cm}: \hspace{0.1cm} ax = b\}
\end{equation*}
where $ax=b$ is a valid inequality for $P$, meaning that it does not cut off any point of polyhedron $P$.
\end{definition}
\begin{definition}\label{facet}
Let $P$ be a non-empty polyhedron in $\mathbb{R}^n$. Face $f$ of polyhedron $P$ is called a facet if the following holds:
\begin{enumerate}
\item[(a)] $f$ is non-empty
\item[(b)] $f$ is not the polyhedron $P$ itself 
\item[(c)] dimension of $f$ is $n-1$.
\end{enumerate}
\end{definition}
From Definition \ref{facet}, it is implied that a non-empty face whose dimension is equal to the dimension of the polyhedron minus one is a facet.
\begin{figure}
    \centering
    \includegraphics[width=0.5\textwidth, height=0.17\textheight]{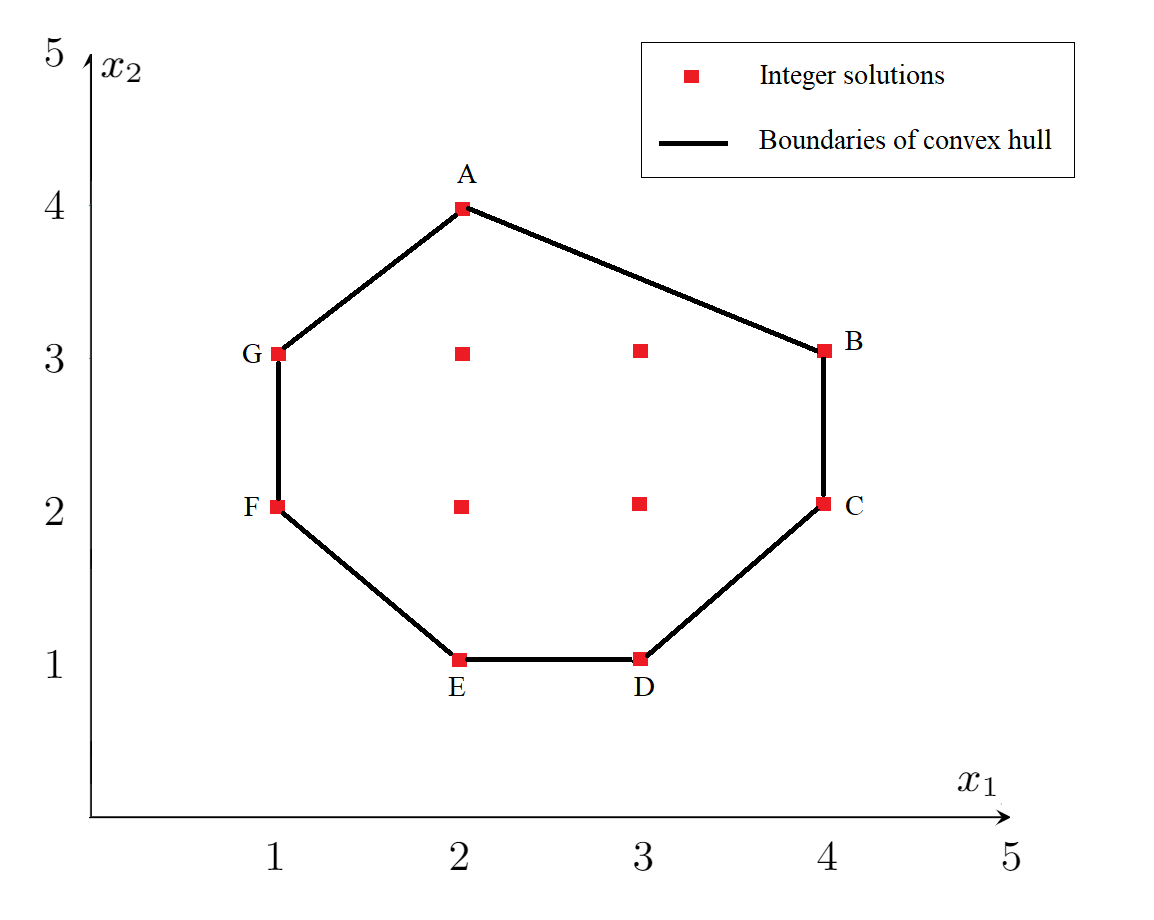}
    \caption{Convex hull, faces, and facets of a set of integer solutions}
    \label{conhull example}
\end{figure}
Figure \ref{conhull example} shows a set of integer solutions with their convex hull in $\mathbb{R}^2$. Let polyhedron $P$ be the convex hull of integer solutions. There are two types of non-empty faces for $P$. First, corner points A, B, C, D, E, F, G are faces of the polyhedron with dimension zero. Second, line segments AB, BC, CD, DE, EF, FG are faces with dimension 1 and therefore are facets. Thus, in the context of linear integer programming, facets are hyper-planes defining a part of the convex hull of integer solutions. If one obtains all the facets of a convex hull, the IP reduces to linear programming and all integrality constraints can be relaxed.

\subsection{Knapsack Inequalities}
Knapsack constraint is a linear constraint of form $\sum_{j \in N} a_jx_j \le b$ where ${a} \in \mathbb{R}_+^N$, ${x} \in \{0,1\}^N$, and $b$ is a positive integer. The polyhedron
\begin{equation}\label{knapsack}
\text{conv} \big\{ x \in \{0,1\}^n \hspace{0.1cm} : \hspace{0.1cm} \sum_{j=1}^n a_jx_j \le b \big\}
\end{equation}
is called a knapsack polytope \cite{wolsey1975faces}. There are many papers on valid inequalities for knapsack polytopes. In particular, minimal cover inequalities (CI) are from the first group of valid inequalities for knapsack polytope which are constructed based on minimal covers. 

A minimal cover $C$ is a set of indices such that $\sum_{j \in C} a_j > b$, and if we remove anyone of these indices, the mentioned inequality is violated, i.e., $\sum_{j \in C \texttt{\textbackslash} \{i\}} a_j \le b, \hspace{0.2cm} \forall i \in C$. For any minimal cover $C$, its minimal cover inequality (CI) can be written as follows \cite{wolsey1975faces}:
\begin{equation} \label{min cover ineq}
\sum_{j \in C} x_j \le |C|-1
\end{equation}
Knapsack CI can be rather
weak. However, they can be strengthened using lifting. Lifting is the process of constructing stronger cuts for an IP (or Mixed IP in general) so that we take a valid inequality and add variables to it such that the generated valid inequality is stronger and even facet-defining in some cases. Balas \citep{balas1975facets} proposed a sequential up-lifting procedure to generate lifted cover inequalities (LCI). They have shown 
such inequalities are stronger than CIs and can even provide facet defining inequalities. However, a key step is to find a minimal CI in this approach, and then the LCI is obtained accordingly.

Recently, \cite{letchford2019lifted} proposed a new lifting procedure that can be applied for cover $C$ even if it is not minimal. The authors shown that the resulting procedure can yield facet-defining inequalities that cannot be obtained by standard lifting procedures, e.g., those of \cite{balas1975facets}. Let set $C$ be a minimal cover. 
Also, let $\overline{a}$ be a positive number such that $\sum_{j \in C} \min \{a_j, \overline{a}\} = b$ where $\overline{a}$ is a unique number which can be calculated using the algorithm proposed in \cite{letchford2019lifted} in $O(|C|)$. Given $\overline{a}$, cover $C$ can be divided in $(C^+,C^-)$ such that $C^- = \{j \in C : a_j \le \overline{a}\}$ and $C^+ = C \texttt{\textbackslash} C^-$. The LCI for any cover $C$, not necessarily minimal, is obtained as follows:

\begin{theorem}\label{theorem 1}\cite{letchford2019lifted}
For all $j \in C$, let $a^-_j = \min \{a_j, \overline{a}\}.$ For $r= 1, \dots, |C|$, let $S^-(r)$ be the sum of the $r$ largest $a^-_j$ values. (Note that $S^-(|C|) = b$.)
Also let $S^-(0)=0$. Finally, given any $k \in N\texttt{\textbackslash} C^-$, let $\gamma_j$ be the largest integer such that $S^-(\gamma_j)\le a_k \le S^-(\gamma_j+1)$. Then the inequality
\begin{equation} \label{new procedure min cover ineq}
\sum_{j \in C^{-}} x_j + \sum_{j \in N\texttt{\textbackslash} C^-} \gamma_j x_j \le |C|-1
\end{equation}
is valid for the knapsack polytope, and it is at least as strong as Balas' lifting procedure.
\end{theorem}

\subsection{Solution Method for The Attack problem}
In what follows, we first show that some integrality constraints of the attack problem can be relaxed while providing a tight relaxation for the problem, and then we will identify constraints for which we can generate LCIs.

Let 
\begin{equation}
F = \text{conv} \Big\{ v^A \in \{0,1\}^{|\mathcal{N}||\mathcal{C}| }, y^A\in \{0,1\}^{|\mathcal{N}|^2}, b\in \{0,1\}^{|\mathcal{C}| }, \alpha \in \mathbb{Z}_+ : \eqref{s1 c1}-\eqref{s1 c9} \Big\}
\end{equation}
be the feasible region for the attack problem. Also, let 
\begin{equation}
F_r = \text{conv} \Big\{ v^A \in \{0,1\}^{|\mathcal{N}||\mathcal{C}| }, y^A\in \{0,1\}^{|\mathcal{N}|^2} : \eqref{s1 c1}-\eqref{s1 c9}, 0 \le b \le 1, \alpha \in \mathbb{R}_+ \Big\}
\end{equation}
be a relaxation of $F$ where integrality restrictions of variable ${\alpha}$ and variable vector $\mathbf{b}$ are relaxed. Considering objective function \eqref{s1 obj func} as $\max_{\mathbb{E}^{A}} f^A(\mathbb{E}^A)$, let 
$$P^A(F) = \{\max_{\mathbb{E}^{A}} f^A(\mathbb{E}^A) : \mathbb{E}^A \in F \}$$
and
$$P^A(F_r) = \{\max_{\mathbb{E}^{A}} f^A(\mathbb{E}^A) : \mathbb{E}^A \in F_r \}$$
be the attack problems considering feasible regions $F$ and $F_r$, respectively. The following proposition shows the equivalence of the two problems.
\begin{proposition} \label{prop equivalence after relaxation}
Problems $P^A(F)$ and $P^A(F_r)$ are equivalent. 
\end{proposition}
\proof{Proof.} Consider problem $P^A(F_r)$. By contradiction, assume variable $\alpha^A$ is not integer and thus $\alpha^{A} = n + \xi$ so that for variables $n$ and $\xi$ we have $n \in \mathbb{Z}^+$ and $ \xi \in (0,1)$. First, constraint \eqref{s1 c2} can be written as  
\begin{equation}\label{relax eq}
\sum_{i \in \mathcal{N}} v^A_{ic} \le \sum_{i \in \mathcal{N}} v^A_{ic'} = \max_{c \in \mathcal{C}} \{ \sum_{i \in \mathcal{N}} v^A_{ic}\} \le \alpha^A = n+\xi
\end{equation}
where $c' \in \mathcal{C}$ is the index corresponding to the maximum value of $\sum_{i \in \mathcal{N}} v^A_{ic}, \forall c$. Since $\sum_{i \in \mathcal{N}} v^A_{ic'}$ is integer as it is the summation of several binary variables, and as a result \eqref{relax eq} is not binding for $c'$, which implies \eqref{relax eq} is not binding for other $c \in \mathcal{C} \backslash \{c'\}$. On the other hand, given at optimality we have $\alpha^{A^*} = n^*+\xi^*$, the optimal value for the objective function \eqref{s1 obj func} is $ f^A(\alpha^{A^*}) = H^* - n^*-\xi^*$ where $H^* = -(|\mathcal{N}|-\sum_{i \in \mathcal{N}} \sum_{c \in \mathcal{C}} v^{A^*}_{ic}) + \sum_{c \in \mathcal{C}} b^{A^*}_c$. Due to the direction of the objective function, an arbitrary solution $\alpha^A_k=n^*+\frac{\xi^*}{k}$ for $k\ge 2$ provides a better value for the objective function since $f^A(\alpha^{A^*})= H^*-n^*-\xi^* < H^*-n^*-\frac{\xi^*}{k}=f^A(\alpha^{A}_k)$. Also, since ${\alpha^A}$ appears only in constraint \eqref{s1 c2},  $\alpha_k^A$ is a feasible solution of $P^A(F_r)$, and the right-hand side of \eqref{s1 c2} becomes $n^*+\frac{\xi^*}{k}$, implying constraint $c'$ is still not binding. Therefore, this is a contradiction implying $\alpha^{A^*}$ is not optimal and $\alpha^{A}_k$ is a better solution. Also as $\frac{\xi}{k}$ converges to zero, constraint \eqref{s1 c2} becomes binding for $c'$, and an upper bound for the value of the objective function is obtained since $\max_k\{H^*-n^*-\frac{\xi^*}{k}\}=H^*-n^*$. Therefore, the optimal value of $\alpha^{A^*}$ in the attack problem is $n^*$ which is integer. A similar argument can be written for constraint \eqref{s1 c3} for $\forall c \in \mathcal{C}$, implying the optimal value of $b$ is also integer. Therefore, by contradiction it is proved that variables $\alpha^A$ and $b$ of the relaxation problem $P^A(F_r)$ take integer values at optimality. Thus, $P^A(F)$ and $P^A(F_r)$ are equivalent. \Halmos

Proposition \ref{prop equivalence after relaxation} shows that one can solve problem $P^A(F_r)$ in the attack stage, which is an equivalent but an easier problem to solve due to relaxing several integrality constraints.

Furthermore, we generate several valid inequalities for the attack problem to reduce to solve time even more. Letting $x_i = \sum_{c \in \mathcal{C}} v^A_{ic}$, the budget constraint \eqref{s1 c5} can be written as a form of knapsack polytope \eqref{knapsack} for which we can derive LCIs based on \cite{letchford2019lifted}. Further details on the implementation procedure can be found in \cite{letchford2019lifted}. We provide further discussions in our numerical results on the efficiency of such LCIs in solving the attack problem.
\subsection{Solution Method for The Response Problem}
The response problem is an IP whose solve time can grow exponentially when the problem becomes large enough. In this section, we propose a novel procedure to reduce the response problem into a smaller and equivalent problem to find the global optimal solution. 

\begin{figure}
    \centering
    \includegraphics[width=0.69\textwidth, height=0.37\textheight]{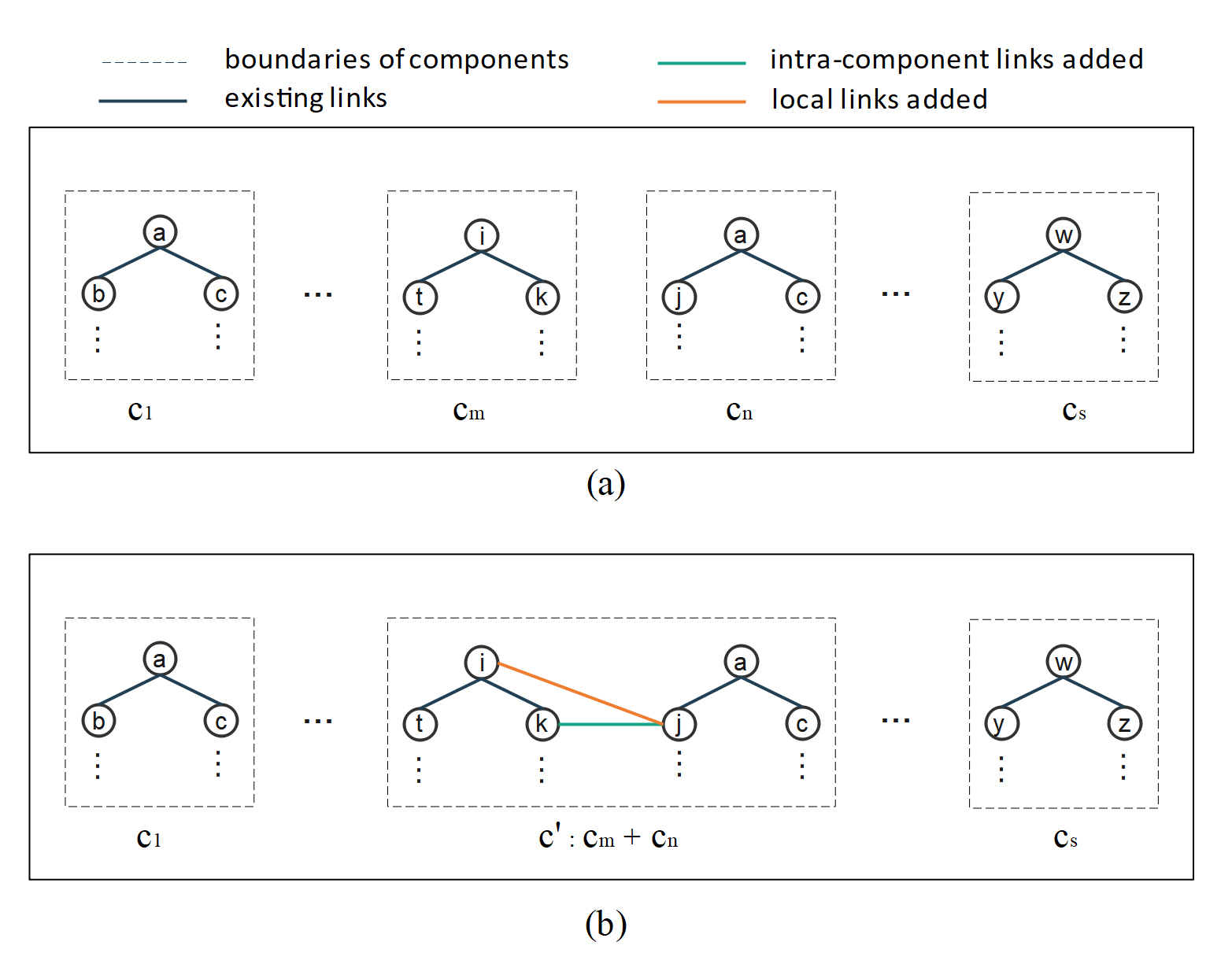}
    \caption{Network topology (a) after attack, (b) after adding intra-component links. The dots represent the remaining of a connected graph in each component.}
    \label{after attack network}
\end{figure}

Figure \ref{after attack network} shows an example of an attacked network consisting of non-empty components $\mathcal{C}= \{c_1, \dots , c_s\}$, where each component consists of a connected subgraph. Let vertices $i$ and $k$ belong to vertices of component $c_m$, ($i,k \in V(c_m)$), and vertex $j$ belongs to vertices of component $c_n$ ($j \in V(c_n))$. There are two candidates for intra-component links as $y^R_{ij}$ and $y^R_{kj}$. Assuming we first add intra-component link $y^R_{kj}$ connecting components $c_m$ and $c_n$, the resulting component becomes a larger component $c' = c_m + c_n$. Therefore, since all vertices of $c_m$ and $c_n$ are now unified in component $c'$, link $y^R_{ij}$ becomes a local link for the new component $c'$, and as shown before, adding a local link would not improve the network connectivity anymore. In other words, from the network connectivity point of view, adding multiple intra-component links between two given components has the same impact on the rupture degree as adding only one intra-component link between the two components. It is because the purpose of link addition is to connect two individual components so that there exists a path between any vertex of the first component and any vertex of the second component. Thus, even though adding multiple intra-component links provides multiple such paths, it does not change the rupture degree which is mainly measured by the number of vertices in the largest component. Therefore, for every two components, all the candidate intra-component links between them have the same value for improving network connectivity but may have different financial values and costs. Thus, to find the optimal solution of the response problem, it is sufficient to consider only the \textit{most cost-efficient intra-component (MCEIC)} links, i.e., the cheapest links, between any two components, one link per two components, to be the candidate links in the response problem. In particular, if we choose to connect two components $i$ and $j$, there is only one candidate MCEIC link between the components which can be added. Doing so, the response problem reduces to a problem determining which MCEIC links should be selected to maximize network resilience given a budget constraint. 

To set up an equivalent problem to the response problem, let $d_{ij}$ be the cost of adding an intra-component link between vertex $i \in V(c_m)$ and vertex $j \in V(c_n)$. Letting 
\begin{equation} \label{matrix cost}
\underline{d}_{mn}=\min_{\substack{i \in V(c_m) \\ j \in V(c_n)}} \{d_{ij}\} , \qquad \forall m,n \in \mathcal{C}^R = \{1, \dots, s\}\end{equation} symmetric matrix $\underline{\mathbf{d}}$ is the min-cost matrix and finds the minimum cost of adding any two components using the MCEIC link between the two components. Also, let $x_{mn} \in \{0,1\}$ be a vector of binary variables indicating 1 if MCEIC link between two components $c_m$ and $c_n$ is added in the reconstructed network. The following proposition shows the response problem reduces to an equivalent problem \eqref{equivalent problem to repair problem}.

\begin{proposition}
The response problem \ref{response problem} is equivalent to 
\begin{subequations}\label{equivalent problem to repair problem}
\begin{align} 
-|X| \hspace{0.1cm} + \hspace{0.1cm} \min_{\mathbf{x}, \mathbf{t}, \alpha^R}  \quad &   \Big( - \alpha^R + {|\mathcal{C}^R| - \sum_{m \in \mathcal{C}^R} \sum_{\substack{n \in \mathcal{C}^R \\ n>m}} x_{mn}} + \epsilon \sum_{m \in \mathcal{C}^R} t^R_m  \Big)  \\
&\sum_{m \in \mathcal{C}^R } t^R_m \ge 1 \label{equ 1}\\
&|V(c_m)| + \sum_{\substack{n \in \mathcal{C}^R \\ n \neq m}} |V(c_n)|x_{mn} \le \alpha^R \le |V(c_m)| + \sum_{\substack{n \in \mathcal{C}^R \\ n \neq m}} |V(c_n)|x_{mn} + M(1-t^R_m)  & \forall m \in \mathcal{C}^R  \label{equ 2}   \\
& \sum_{m \in \mathcal{C}^R} \sum_{\substack{n \in \mathcal{C}^R \\ n>m}} \underline{d}_{mn}x_{mn} \le B^R & \label{equ 3}\\
& x_{mn} = x_{nm} \hspace{8cm} \forall m,n \in \mathcal{C}^R \label{equ 4} \\ & \mathbf{x} \in \{0,1\}^{| \mathcal{C}^R|^2}, \mathbf{t} \in \{0,1\}^{|\mathcal{C}^R|}, \alpha^R \in \mathbb{R}_+ \label{domain aux} 
\end{align}
\end{subequations}
where $|X| = (|\mathcal{N}|-\sum_{i \in \mathcal{N}} \sum_{c \in \mathcal{C}} v^{A*}_{ic})$ is the size of the worst cut set which was identified in the attack problem.
\end{proposition}
\proof{Proof. } 
Constraints \eqref{s2 c1}, \eqref{s2 c6.1}-\eqref{eq4} of response problem \ref{response problem} correspond to the structure of components to ensure only intra-component links can be added in the reconstructed network. It easily follows that these constraints are all implied in Formulation \eqref{equivalent problem to repair problem} using the definition of MCEIC links. Also, to find the number of components, we know each MCEIC link unifies two components and thus reduces the number of components by 1. Therefore, $(|\mathcal{S}| - \sum_{m \in \mathcal{S}} \sum_{\substack{n \in \mathcal{S} \\ n>m}} x_{mn})$ corresponds to the number of components after link addition, which is considered in the objective function, and implies that constraint \eqref{s2 c3} is not necessary anymore. Furthermore, allowing the MCEIC links to be the only candidate links in the attack problem, constraints \eqref{s2 aux2}, \eqref{s2 aux1}, \eqref{budget}, \eqref{s2 c6} of the response problem correspond to \eqref{equ 1}, \eqref{equ 2}, \eqref{equ 3}, \eqref{equ 4}, respectively. In addition, after relaxing restrictions on $\alpha$ (based on Proposition \ref{prop equivalence after relaxation}), constraint \eqref{domain variable} reduces to \eqref{domain aux}. \Halmos

Given the optimal solution of \eqref{equivalent problem to repair problem}, one can find the intra-component links added to the reconstructed network ($\mathbf{y}^R$), the length of the largest component ($\alpha^R$), the number of components, and thus, the optimal value of network resilience under the worst-case attack. It is to be note Formulation \eqref{equivalent problem to repair problem} is less complicated than the response problem \ref{response problem} and has significantly a less number of constraints and integer variables, and thus is expected to be computationally less challenging. Furthermore, we can show that there exists an embedded knapsack constraint in problem \ref{equivalent problem to repair problem} for which we can generate LCIs so that solving formulation \ref{equivalent problem to repair problem} becomes even more time-efficient. To do so, we first present a simple re-indexing procedure to convert two-dimensional matrices into vectors, which can be done in polynomial time.


\subsubsection{Re-indexing Procedure}
In Formulation \ref{equivalent problem to repair problem}, matrix variable $\mathbf{x}$ is a two-dimensional symmetric matrix where $x_{mn}=x_{nm}$ and $x_{nn}=0, \forall m,n \in \{1, \dots, s\}$. Let $|\mathcal{C}^R|=s$. To identify the underlying knapsack constraint from Formulation \ref{equivalent problem to repair problem}, it is sufficient to take into consideration the upper triangular matrix $\mathbf{x}$ which has $s(s-1)/2$ elements, instead of considering $s^2$ elements in matrix $\mathbf{x}$. Then, we store all elements of upper triangular matrix $\mathbf{x}$ in a vector $\hat{x}$.
\begin{figure}
    \centering
    \includegraphics[width=\textwidth]{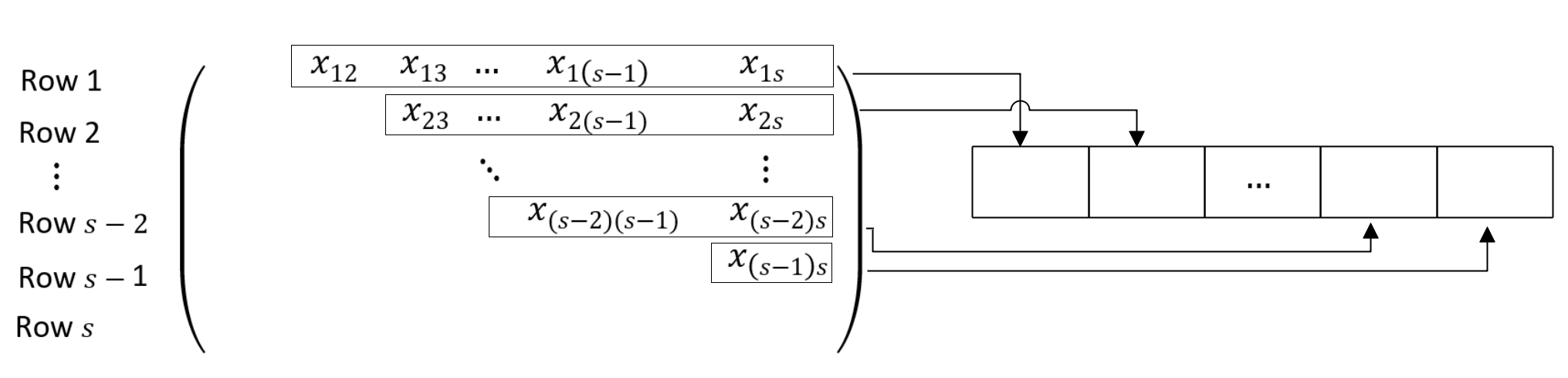}
    \caption{Mapping upper triangular matrix $\mathbf{x}$ to vector $\hat{\mathbf{x}}$.}
    \label{fig: mapping}
\end{figure}
Figure \ref{fig: mapping} shows the re-indexing and mapping from matrix $\mathbf{x}$ to vector $\hat{\mathbf{x}}$. Particularly, we place all rows of matrix $\mathbf{x}$ in vector $\hat{\mathbf{x}}$ where row $m$ of matrix $\mathbf{x}$ starts at position $S_I(m)=1+\sum_{k=1}^{m-1}(s-k)$ of vector $\hat{\mathbf{x}}$. Therefore, considering any element $x_{mn}$ in the upper triangular matrix $\mathbf{x}$  ($m < n$), this element is mapped to the $\sigma(m,n)$-th element of vector $\hat{\mathbf{x}}$, denoted by $\hat{x}_{\sigma(m,n)}$, where
\begin{align}
\sigma(m,n)= (n-m)+\sum_{k=1}^{m-1}(s-k) \label{re-index position}
\end{align}
Given this re-indexing, for example $\sum_{\substack{n \in \mathcal{C}^R \\ n \neq m}} |V(c_n)|x_{mn}, \forall m \in \mathcal{C}^R$ in \eqref{equ 2} can be re-written as follows
\begin{equation}
\sum_{n=1}^{m-1} |V(c_n)|\hat{x}_{\sigma(n,m)} + \sum_{n=m+1}^{s} |V(c_n)|\hat{x}_{\sigma(m,n)} \qquad \quad \forall m \in \mathcal{C}^R
\end{equation}

Similarly, let vector $\tilde{\mathbf{d}}$ be a cost vector with the same size as vector $\hat{x}$ where $\tilde{\mathbf{d}}=\{\underline{d}_{12},\dots, \underline{d}_{1s}, \dots, \underline{d}_{rs} \dots, \underline{d}_{(s-1)s}\}$, implying that the cost of adding a MCEIC link between components $m$ and $n$ is stored in $\tilde{d}_{z}$ where $z=\sigma(m,n)$. Thus, given the re-indexing procedure, the Formulation \ref{equivalent problem to repair problem} reduces to Formulation \ref{final - equivalent problem to attack problem}.
\begin{subequations}\label{final - equivalent problem to attack problem}
\begin{align} 
-|X| + \min_{\hat{\mathbf{x}}, \mathbf{t}, \alpha^R}  \quad &  - \alpha^R + (|\mathcal{S}| - \sum_{z = 1}^{s(s-1)/2} \hat{x}_{z}) - \epsilon \sum_{m \in \mathcal{S}} t^R_m \\
&\sum_{m=1}^s t^R_m \ge 1 \label{final - equ 1}\\
&|V(c_m)| + \sum_{n=1}^{m-1} |V(c_n)|\hat{x}_{\sigma(n,m)} + \sum_{n=m+1}^{s} |V(c_n)|\hat{x}_{\sigma(m,n)} \le \alpha^R  \nonumber \\ & \le |V(c_m)| + \sum_{n=1}^{m-1} |V(c_n)|\hat{x}_{\sigma(n,m)} + \sum_{n=m+1}^{|\mathcal{S}|} |V(c_n)|\hat{x}_{\sigma(m,n)} + M(1-t^R_c)  & \forall m \in \mathcal{C}^R   \label{final - equ 2}   \\
& \sum_{z=1}^{s(s-1)/2}  \underline{d}_{z}\hat{x}_{z} \le B^R & \label{final - equ 3}\\
& \hat{\mathbf{x}} \in \{0,1\}^{s(s-1)/2}, \mathbf{t} \in \{0,1\}^{s}, \alpha^R \in \mathbb{R}_+ \label{final - domain aux}
\end{align}
\end{subequations}
Formulation \ref{final - equivalent problem to attack problem} is a smaller problem than Formulation \ref{equivalent problem to repair problem} due to reducing the number of variables and constraints by the re-indexing procedure. Furthermore, constraint \eqref{final - equ 3} corresponds to a knapsack constraint for which we derive LCIs based on \cite{letchford2019lifted}, to reduce the computational challenges. As a result, instead of the response problem \ref{response problem}, one can solve Formulation \ref{final - equivalent problem to attack problem} which is time-efficient while providing the same global optimal solution. 

\subsection{Extension to Other Types of Attacks} \label{section other attacks}
There are several types of attacks studied in the literature of network resilience and robustness. A common approach is to identify the worst-case attack and optimize the network under such attacks as it was shown  that cascading failures are mainly due to targeted attacks \citep{ruj2014analyzing}. 

However, there are other types of attacks or failures in various networks. For instance, in several networks, such as power system networks or supply chain networks, the random failure of nodes plays an important role in network non-resilience due to the probabilistic age function of equipment, natural disasters, etc. \citep{sadrabadi2021resilient,ji2016improving}. We refer to such failures as \textit{random attacks}.
Furthermore, in some cases, the attackers desire to disable a set of certain nodes, not necessarily the most critical nodes, in a network due to geographical location, importance, political reasons, etc. We refer to such attacks as \textit{designated attacks}. Finally, in some networks, the nodes can be divided into some subgroups where the attacker is interested in attacking certain subgroups based on their properties. For example, the nodes of a power system network can be divided into generation and load nodes, where if the generation nodes are disabled, the load nodes cannot be supplied automatically. Thus, the attacker might be interested in disabling the generation nodes to cut the power supply. We refer to this type of attack as \textit{distributed attacks}.
\subsubsection{Random and Designated Attacks}
In the case of having random, one needs to determine the set of nodes that are most likely to be attacked before optimizing the network resilience. Such analysis can be done using probabilistic evaluation based on the failure distributions of nodes, Monte Carlo simulation, Markov Chain models, etc. Similarly, in case of having designated attacks, the nodes that are subject to attack are known in advance. 

Let $X^{\text{n}}$ be the cut set corresponding to the random or designated attacks. If $G(V-X^{\text{n}})$ is a connected graph, there is a path between any two nodes of the attacked network and all the remaining nodes are able to communicate. Thus, the network resilience is maximized. However, if $G(V-X^{\text{n}})$ is not connected, once can apply the proposed framework to maximize network resilience. Since these cut sets are known or estimated in advance, the network resilience is to be maximized under such cut sets, as opposed to the worst cut set. Thus, the attack problem is not needed to be solved anymore. As a result, our proposed two-stage formulation reduces to only the response stage where the network resilience is aimed to be maximized under a pre-identified cut set. By letting $X=X^{\text{rand}}$ or $X=X^{\text{d}}$ in the response problem, our proposed approach is able to optimize network resilience for random or designated attacks.

\subsubsection{Distributed Attacks} \label{distributed attacks}
In the case of having distributed attacks, one needs to classify each node of a network in subsets of nodes considering a certain property of interest. Without loss of generality, we divide the nodes in two subsets and it can be extended to more subsets. 
Consider a network with a set of nodes denoted by $\mathcal{N}$. Let $\mathcal{N}_I$ and $\mathcal{N}_J$ be two subsets of $\mathcal{N}$ such that $\mathcal{N}_I \cap \mathcal{N}_J = \emptyset$ and $\mathcal{N}_I \cup \mathcal{N}_J = \mathcal{N}$. Supposing the attacker is interested in only disabling the nodes of subset $\mathcal{N}_I$. Thus, Constraint \ref{s1 c1} is re-written as follows:
\begin{subequations}
\begin{align}
& \sum_{c \in \mathcal{C}} v^A_{ic} \le 1 & \forall i \in \mathcal{N}_I \label{dist s1 c1} \\
& \sum_{c \in \mathcal{C}} v^A_{ic} = 1 & \forall i \in \mathcal{N}_J \label{dist s2 c1}
\end{align}
\end{subequations}
where Constraint \eqref{dist s1 c1} allows the possibility of disabling each vertex $i \in \mathcal{N}_I$ by the worst-case attack. Otherwise, the vertex remains active and is at most in one component. Also, Constraint \eqref{dist s2 c1} enforces that all vertices of subset $\mathcal{N}_J$ remain active after the attack and belong to exactly one component after the worst-case attack. This reformulation can also be used for the type of attacks mentioned in \cite{abbas2017improving} where a subset of nodes always remain intact and cannot be disabled. The rest of the proposed two-stage model remains the same as before and one can use the proposed solution method to address the network resilience and robustness under such distributed attacks.

\section{Numerical Results}\label{section; results}

In this section, we first provide a illustrative numerical example to show the benefits of our methodology in terms of network resilience and robustness and the effects of the budget on the link addition strategies. Then, we use larger problems to demonstrate the computational benefits of the proposed solution method. In all numerical examples, we generate random connected graphs and execute the proposed algorithm. For all instances, we let the cost for attacking each node be 1 unit and the total budget for attack be $|\mathcal{N}|/2$ units so that no more than half of the nodes are allowed to be disabled in the attack stage. All instances are performed using the MIP solver of Gurobi Version 9.1 on a 1.5 GHz×86 laptop with 1GB main memory, and the optimality gap is set to 0 for all instances.

\subsection{Resilience and Robustness Analysis}
\begin{figure}
    \centering
    \includegraphics[width=0.77\textwidth, height=0.30\textheight]{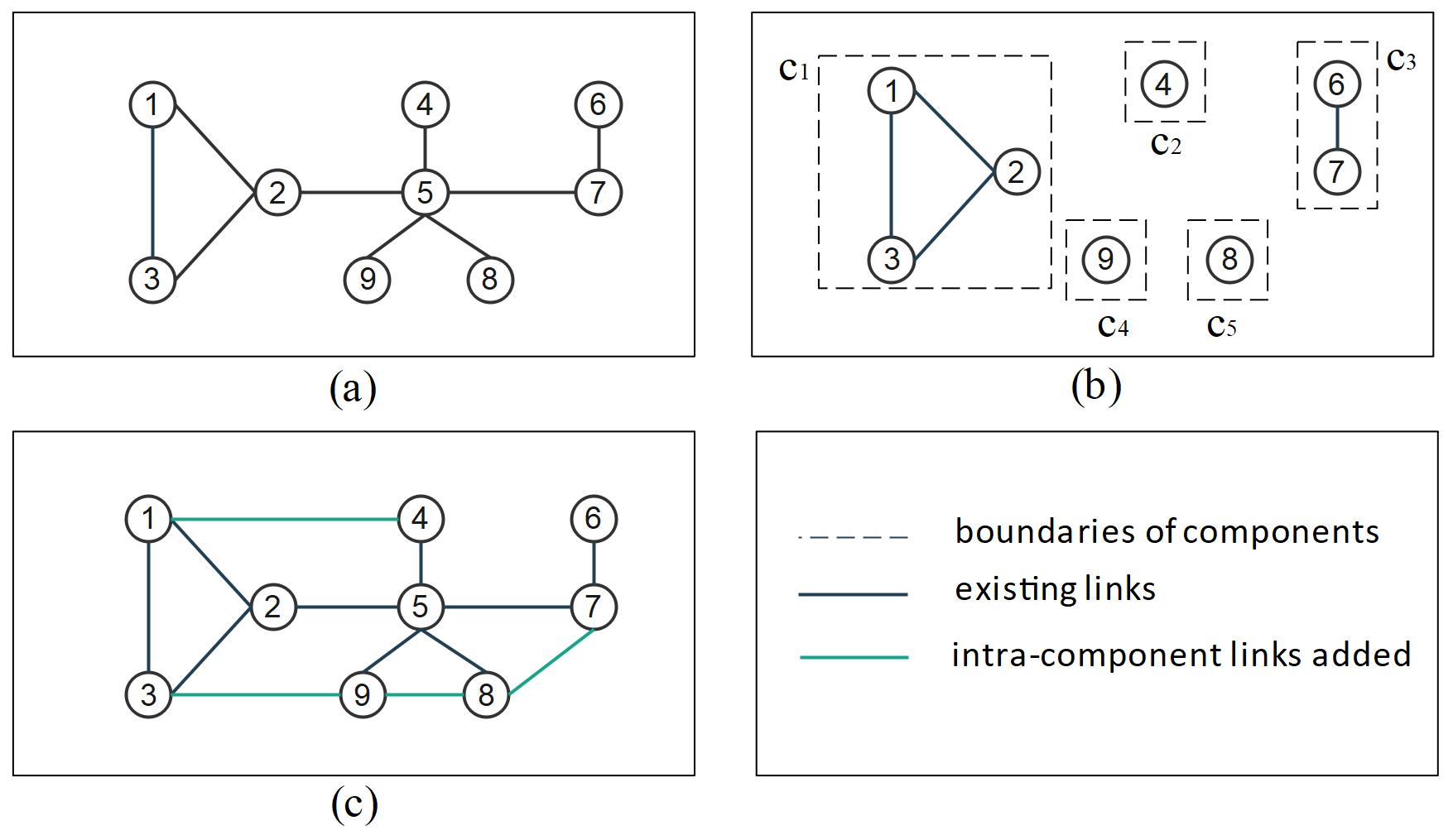}
    \caption{The topology of (a) the initial network, (b) the attacked network (stage I), (c) the reconstructed network (Stage II)}
    \label{fig example}
\end{figure}

In this section, we present a network with 9 vertices and 9 edges, shown in Figure \ref{fig example}.(a), to demonstrate the benefits of the proposed network. Solving the attack problem, the worst cut set is $X=\{5\}$ corresponding to the lowest network resilience ($-r_X(G)$). This implies the connectivity of this network heavily relies on vertex 5. Therefore, the worst-case attack can occur by disabling vertex 5. Figure \ref{fig example}.(b) shows the structure of the network under the worst-case attack by removing vertex 5. Since there exist 5 components in Figure \ref{fig example}.(b), one can interpret that the maximum network resilience can be obtained by adding 4 intra-component links, if the link addition budget allows. It is because once all 5 components are connected together using 4 intra-component links, the resulting network consists of only one component. Figure \ref{fig example}.(c) shows the reconstructed network under the worst-case attack which consists of the original network and the new intra-component links added to the network. For now, we do not consider the budget constraint for link addition. The green links are those added to maximize the network resilience under the worst-case attack to reduce the dependencies on vertex 5 so that the remaining vertices are able to properly communicate in the absence of vertex 5. 

In Table \ref{Table: cut sets}, we manually choose a number of cut sets $X_i$s, not necessarily the worst cut set, and measure the network resilience
of the initial network, denoted by $-r_X(G)$ and that of the reconstructed network, denoted by $-r^R_X(G)$. As observed, in the initial network, $X_1=\{5\}$ is the worst cut set corresponding to the lowest value of the network resilience which shows the network is more vulnerable and can seriously be damaged in the absence of vertex 5. For other cut sets, the network is more resilient as $-r_{X_1}(G)=-1 \le -r_{X_i}(G),  \forall i \neq 1$, implying
either the network is less dependent on the vertices of other cut sets or there is more work required to damage the network by disabling more vertices of a cut set.

\begin{table}[!t]
\centering
\footnotesize
\caption{Resilience and Robustness analysis for different cut sets of the initial and the reconstructed network}
\label{Table: cut sets}
\scalebox{0.9}{
\begin{tabular}{lcccc}
\hline
& \multicolumn{2}{c}{Initial Network}&\multicolumn{2}{c}{Reconstructed Network}   \\ \hline
 Cut set & \# of components & $\text{NR}: -r_X(G)$ & \# of components & $\text{NR}: -r^R_X(G)$ \\
\hline
$X_1=\{5\}$ & 5 & \textbf{ -1} $(r^*)$ & 1 & 8 \\ 
$X_2=\{2\}$ & 3 & 5 & 1 & 8 \\
$X_3=\{5,7\}$ & 4 & 0 & 2 & 6 \\
$X_5=\{2,5,7\}$ & 5 & 0 & 2 & 6 \\
$X_4=\{2,5,8\}$ & 4 & 1 & 2 & 5 \\
$X_5=\{2,5,7,9\}$ & 4 & 2 & 3 & 4 \\
$X_6=\{1,2,5,9\}$ & 5 & 2 & 3 & 4 \\
$X_7=\{1,2,5,8\}$ & 4 & 2 & 3 & 3 \\
$X_8=\{1,2,5,7,9\}$ & 4 & 2 & 4 & \textbf{  2} $(r^{R*})$  \\
$X_9=\{3,2,5,7,9\}$ & 4 & 2 & 3 & 4  \\
\hline
\multicolumn{3}{l}{NR: Network Resilience} \\
\multicolumn{3}{l}{*: optimal objective value of the attack problem} \\
\hline
\end{tabular} }
\end{table}
From Table \ref{Table: cut sets}, there are three main observations on the reconstructed network, which are obtained after immunizing the network against the worst-case attack. First, the network resilience for $X_1=\{5\}$ significantly improves from $-1$ to $8$ after we add links meaning that the functionality of the reconstructed network does not heavily depend on only vertex 5. Thus, if vertex $5$ is damaged, the rest of the network is able to maintain its proper functioning and communication with other parts of the network. Therefore, the reconstructed network is more resilient in the absence of the worst-case attack. Second, after adding links to the initial network, the network dependencies on only vertex 5 decrease. Thus, the new worst cut set in the reconstructed network may consist of other vertices along with vertex 5. We call this \textit{dynamic} worst cut set since it is obtained after the link addition strategies are known. In particular, the cut set $X_8=\{1,2,5,7,9\}$ becomes the dynamic worst cut set for the reconstructed network. This means that after links are added to the network, the most severe damage to the reconstructed network occurs only when all vertices 1, 2, 5, 7, and 9 are simultaneously disabled. We refer to this as the network robustness measured by the size of the worst cut set. Thus, by adding links, the network robustness increases from $|X_1|=1$ to  $|X_8|=5$ implying that disabling all these five vertices simultaneously requires more budget and more work done to damage the network, which might not be a viable option for the attacker due to the budget constraint considered in the attack problem. 
Third, assuming the attacker can dynamically respond to our link addition actions and has enough budget to disables all vertices of cut set $X_8$, we observe that the network resilience for the reconstructed network under the new (dynamic) worst-case attack is more than that of the initial network, i.e., $ r^* =-1 < r^{R*}=2$. This shows that the reconstructed network not only is significantly more resilient against the worst cut set of the initial network, i.e., $X_1$, but also it becomes more resilient against the dynamic worst cut set of the reconstructed network, i.e., $X_8$, that may be identified if the attackers can dynamically change their decisions based on our link addition strategy. The mentioned observations demonstrate how the proposed approach leads to a reconstructed network that is both more resilient and robust against the initial and dynamic worst-case attacks.

%
%
\subsection{Effects of the link addition budget}
In this section, we show the effects of the budget considered for link addition on the network resilience and robustness. To do so, we first find the matrix $\underline{\mathbf{d}}$ obtained from \eqref{matrix cost} to find all candidate MCEIC links of the attacked network considered in Figure \ref{fig example}.b. The cost of adding MCEICs link and the two ends of such links, presented in parenthesis, are shown in \eqref{dmin example}.
\begin{equation} \label{dmin example} 
\allowdisplaybreaks
\underline{\mathbf{d}} = \begin{blockarray}{cccccc}
c_1 & c_2 & c_3 & c_4 & c_5 \\
\begin{block}{(ccccc)c}
  - & 1.5 (1\textendash4) & 2.2 (2\textendash7) & 1.5 (3\textendash9) & 1.6 (2\textendash8) & c_1 \\
   & - & 1.8 (4\textendash6) & 2.2 (4\textendash9) & 2.2 (4\textendash8) & c_2 \\
   &  & - & 2.3 (7\textendash9) & 1.3 (7\textendash8) & c_3 \\
   &  &  & - & 1.0 (8\textendash9) & c_4 \\
   &  & &  & - & c_5 \\
\end{block}
\end{blockarray}
 \nonumber
\end{equation}
%
\begin{figure}
\centering
\medskip
\begin{tikzpicture}[font=\small]
\pgfplotsset{width=0.49\textwidth, height=0.25\textheight}
\begin{axis}[legend columns=1, legend image post style={scale=1},legend style={/tikz/every even column/.append style={column sep=0.2cm},{font=\scriptsize}, at={(0.0,0.761)}, anchor=south west},
xmin=0,
    xmax=6,
    xtick={1,2,3,4,5,6},
    xlabel= $B^R$,
    ymin=-2,
    ymax=12,
    ytick={-2,0,2,4,6,8,10,12},
    ylabel={Network Resilience}]
\addplot [no marks, black,thick, name path=conventional approach] plot coordinates { (0,-1 ) (1, -1) (1, 0) (1.5, 0) (1.5, 1)  (2.2,1) (2.2,2) (2.5,2) (2.5,3) (2.9,3)(2.9,4) (3.8,4) (3.8,6) (5.3,6) (5.3,8) (6,8)
};

\addplot [no marks, black, name path=h0] plot coordinates { (0,-2) (0,12)
};
\addplot [no marks, gray!5, name path=h1] plot coordinates { (1.5,-2) (1.5,12)
};
\addplot [no marks, gray!10, name path=h2] plot coordinates { (2.5,-2) (2.5,12)
};
\addplot [no marks, gray!20, name path=h3] plot coordinates { (3.8,-2) (3.8,12)
};

\addplot [no marks, gray!25, name path=h4] plot coordinates { (5.3,-2) (5.3,12)
};

\addplot [no marks, black, name path=h5] plot coordinates { (6,-2) (6,12)
};

\addplot [no marks,dashed, black,thick, name path=h5] plot coordinates { (106,-2) (106,12)
};

\addplot[gray!5] fill between[ 
    of = h0 and h1];
\addplot[gray!12] fill between[ 
    of = h1 and h2];    
\addplot[gray!19] fill between[ 
    of = h2 and h3];
\addplot[gray!26] fill between[ 
    of = h3 and h4];
\addplot[gray!33] fill between[ 
    of = h4 and h5];
%
%
%
%
\legend{ $-r^R_{X_1}(G)$,,,,,,,$|X^*|$ }
\end{axis}
\begin{axis}[legend image post style={scale=0.6},legend style={{font=\tiny}, at={(0.2081,0.0)}, anchor=south west},
   xmin=0,
    xmax=6,
axis x line=none,
    axis y line*=right,
    ymin=-4,
    ymax=6,
    ytick={1,2,3,4,5,6},
    ylabel={Network Robustness}
  ] 
\addplot [no marks, dashed, black,thick, name path=conventional approach] plot coordinates { (0,1 ) (1.5, 1) (1.5, 2)  (2.5,2) (2.5,3)(3.8,3)(3.8,4)(5.3,4) (5.3,5) (6,5) };
\end{axis} 
\end{tikzpicture}
\caption[Short caption]{Network resilience (shown by the left vertical axis) and robustness (shown by the right vertical axis) for different values of link addition budgets. Shaded areas correspond to various ranges of the budget with a constant number of links added where the brightest area has 0 links, the darkest area has 4 links, and the increment is 1 link.}
    \label{fig effect of budget}
\end{figure}
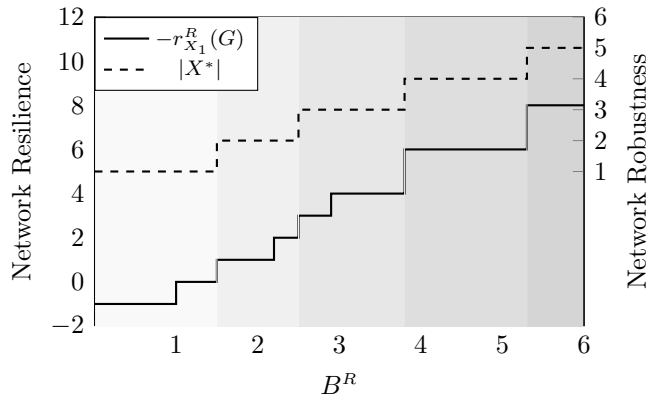
%
After the worst attack $X_1=\{5\}$ is identified, one can maximize the network resilience by link addition within a budget so that the network maintains a maximum functionality level in the absence of node $5$. 
However, when the budget is limited, one cannot connect all the remaining components. Thus, the best combination of the intra-component links should be selected to maximize the network resilience within the budget. In addition, as we showed, the link addition reduces the dependencies on the worst cut set $X_1=\{5\}$, and includes more vertices in the new cut set for the reconstructed network. Thus, we expect to have larger cut sets (higher robustness) when more intra-component links are added to the network. 

Figure \ref{fig effect of budget} demonstrates the effects of link addition budget on the reconstructed network under the initial worst-case attack $X_1=\{5\}$. The network resilience (measured by $-r_{X_1}^R(G)$ and shown by the left vertical axis) and the network robustness (measure by $|X^*|$ and shown by the right vertical axis) for various values of the budget assigned for the link addition. In the figure, shaded areas correspond to various ranges of the budget within which a constant number of links is added, where the brightest area on the left indicates 0 links are added in the reconstructed network, the darkest area indicates 4 links are added in the reconstructed network, and the increment is 1 link. For example, $ 2.5 \le B^R < 3.8$ corresponds to the third shaded area where two MCEIC links are added to the reconstructed network. From the figure, it is observed that once the budget is less than 1, no MCEIC link is added in the repair stage since the budget is not enough. Thus, the reconstructed network is equivalent to the initial network with $-r^R_{X_1}(G)=-1$. However, as the budget increases, the best MCEIC link candidates given different levels of the budget are added to the network which corresponds to a higher value for $-r^R_{X_1}(G)$, which in turn enhances the network resilience under the initial worst-case attack. For example, when the budget is $1.3$ units, the only MCEIC link that can be added is link $7$-$8$, which unifies components $c_3$ and $c_5$ (decreases the number of components by 1) without changing the largest component $c_1$. Thus, $-r^R_{X_1}(G)=0$. However, when the budget increases to 1.5, the response problem allows to add MCEIC link $1$-$4$, instead of link $7$-$8$, since the new link connects components $c_1$ and $c_3$ (drops the number of components by 1) while increasing the length of the largest component $c_1$ by 1. Thus, $-r^R_{X_1}(G)=1$. By further changing the budget assigned for link addition, the optimal topology of the reconstructed network changes with objective to increase the network resilience. Finally, once there is sufficient budget (i.e., any budget greater than 5.3 units), 4 MCEIC links are added to the network to unify all components resulting in maximum resilience, which corresponds to the solution of the previous part where we ignored the budget constraint.

The secondary axis of Figure 
\ref{fig effect of budget} evaluates the network robustness $|X^*|$ of the reconstructed network for various levels of the budget. To evaluate the network robustness, we dynamically find the new worst cut sets of the reconstructed network, as the dependencies may change once new links are added. For the initial network, $|X^*|=1$, meaning that by removing one vertex the network experiences the more severe damage. However, as the budget increases the reconstructed network becomes more robust so that one may need to disable more number of vertices to severely damage the network, which might not be an affordable option for the attacker. Ultimately, once there is a sufficient budget, i.e., $B^R \ge 5.3$, the dynamic worst cut set includes five vertices implying that damaging the network is more difficult since all five vertices should be simultaneously disabled. The mentioned analyses provide an intuitive picture for the managers to determine the desired level of robustness and resiliency that they would like to consider within a financial budget.

\subsection{Computational and Robustness Analysis}
In this section, we further evaluate the performance of the proposed framework, mostly focusing on the computational results. We randomly generate 14 instances to execute the proposed framework\footnote{ the data for the random instances are available at  \url{https://github.com/miladdf94/Network-resilience/blob/main/Adj.\%20matrix\%20-\%20instances.txt}}. In our implementation, we relax the link addition budget constraint (in the response problem) to let the proposed framework maximizes the network resilience without financial limitations.

Table \ref{Table: cut sets 1} and \ref{Table: cut sets 2}, respectively, show the resilience and robustness results as well as computational results for 14 randomly-generated instances. From Table \ref{Table: cut sets 1}, it can be observed that:
\begin{itemize}
    \item For all instances, the network robustness noticeably increases, where $|X^*|$ and $|X^{'*}|$ measure the robustness of the initial and the reconstructed networks. In particular, for instance \#1, $|X^*|$ and $|X^{'*}|$ are 3 and 6, respectively, implying that the most severe damage in the initial network occurs when 3 vertices are attacked simultaneously, while to most severely attack the reconstructed network, 6 vertices are required to be disabled simultaneously.
    \item In all instances, the reconstructed network is significantly more resilient against the initial worst-case attack (denoted by $-r^R_{X^*}(G)$) compared to the initial network (denoted by $-r_{X^*}(G)$). In addition, if the attacker has the ability to change its decisions based on the link addition strategy in the response stage, the new worst-case attack is denoted by $X^{'*}$, which is impacted by the new links added to the initial network. It can be observed that in such a case, the network resilience in the reconstructed network is still significantly higher than that of the initial network, i.e., $ -r_{X^*}(G) < -r^R_{X^{'*}}(G)$ for all instances, demonstrating the benefits of the proposed framework both for static and dynamic worst-case attacks.
\end{itemize}

Also, Table \ref{Table: cut sets 2} presents the computational results of the mentioned instances. We first solve the proposed two-stage model directly in Gurobi, and then implement our proposed solution method and solve the two-stage model in Gurobi. In Table  \ref{Table: cut sets 2}, column "solver cut" corresponds to the number of cuts generated in Gurobi, and "add. user cuts" indicated the number of additional LCI cuts \citep{letchford2019lifted} that we manually added to the solver cuts to solve the attack problem. Also, "nodes" corresponds to the number of nodes (branches) that are explored in the branch-and-bound (B\&B) algorithm to find the optimal solution. For the response problem, we did not add LCI cuts since we relaxed the link addition budget constraint to more clearly demonstrate the performance of the proposed method. 

From Table \ref{Table: cut sets 2}, it can be observed that the proposed solution method significantly reduces the solve time by exploiting the problem structure. For the attack problem, the proposed approach solves the attack problem faster and more efficiently. This is due to the efficient relaxation in the proposed approach as well as the LCIs added to the attack problem. On the other hand, in the response stage, the proposed approach corresponds to a much simpler, but equivalent, problem that can be solved faster than the original response problem. In particular, our proposed solution method requires exploring 1 node in the B\&B algorithm to find the optimal solution. However, without implementing the proposed algorithm, the B\&B algorithm exhaustively explores a large number of nodes to find the optimal solution. The mentioned observations validate the benefits of the proposed solution method from the resilience and robustness point of view as well as computational aspects.
\begin{table}[!t]
\centering
\footnotesize
\caption{Resilience and Robustness analysis of initial and the reconstructed network for various instances}
\label{Table: cut sets 1} 
\scalebox{1}{
\begin{tabular}{c|cc|cc|cc|ccc}
\hline
& \# of & \# of & Used & \# of MCEIC  & \multicolumn{2}{c|}{Initial network} & \multicolumn{3}{c}{reconstructed network} \\ \cline{6-7} \cline{7-10}
Instance & nodes & edges & budget & links & $|X^*|$ & $-r_{X^*}(G)$ & $-r_{X^*}^R(G)$& $|X^{'*}|$ &  $-r_{X^{'*}}^R(G)$ \\ \hline 
1&11&15&5.55&4&3&0&10&6&6 \\   
2&11&20&6.95&4&5&2&10&6&6 \\    
3&12&20&6.70&5&5&1&11&6&6 \\ 
4&12&30&4.50&3&6&4&11&7&7 \\ 
5&15&19&14.05&8&5&-2&14&6&6 \\
6&15&26&12.05&6&7&2&14&8&8\\
7&16&20&9.05&7&5&-1&15&7&7 \\
8&16&32&7.50&5&7&3&15&8&8 \\
9&17&21&11.40&8&6&-1&16&7&2 \\
10&17&30&9.65&7&7&1&16&9&9 \\
11&18&19&12.25&8&5&-2&17&7&7 \\
12&18&25&13.75&8&6&-1&17&8&8 \\
13&21&25&14.10&10&7&-2&21&8&8 \\
14&22&27&11.75&10&7&-2&21&9&9 \\
\hline
\end{tabular} }
\end{table}

\begin{table}[!t]
\centering
\footnotesize
\caption{Computational analysis of initial and the reconstructed network for various instances}
\label{Table: cut sets 2}
\scalebox{0.87}{
\hspace{-0.8cm}\begin{tabular}{c|ccc|cc|c|ccc|cc|c}
\hline
& \multicolumn{6}{c|}{without the proposed solution method} & \multicolumn{6}{c}{with the proposed solution method} \\ \cline{2-13} 
& \multicolumn{3}{c|}{Attack problem} & \multicolumn{2}{c|}{Response problem} & Total & \multicolumn{3}{c|}{Attack problem} & \multicolumn{2}{c|}{Response problem}  & Total\\ \cline{2-6} \cline{8-12}
Instance & solver cuts & nodes & time (s) & nodes & time (s) &time (s)& add. user cuts & nodes & time (s) & nodes & time (s) & time (s)\\ \hline
1&120&3723&2.14&8&0.12&2.26&14&786&0.98&1&0.03&1.01 \\
2&129&4665&2.64&156&0.14&2.79&11&937&1.10&1&0.02&1.12 \\
3&143&3613&3.26&206&0.14&3.40&19&704&0.93&1&0.02&0.95\\
4&98&10757&5.69&675&0.28&5.96&7&1589&1.92&1&0.03&1.95\\
5&248&8374&6.80&5093&2.13&8.93&18&1396&1.74&1&0.03&1.77\\
6&728&118584&60.21&12034&7.07&67.29&35&60543&36.10&1&0.03&36.13\\
7&523&40959&27.15&798&4.15&31.31&22&34367&25.19&1&0.03&25.22\\
8&451&281822&177.78&23092&12.06&189.54&	81&	116553&	89.20&	1&	0.02&	89.21\\
9&750&	157823&	91.25&	16798&	11.16&	102.41&	64&	74979&	51.47&	1&	0.01&	51.49\\
10&663&	131102&	103.38&	12489&	9.09&	112.47&	14&	57593&	54.62&	1&	0.01&	54.63\\
11 & 329&	3827887&	1416.52&	78920&	28.98&	1445.50&	103&	1975251&	857.02&	1&	0.01&	857.03\\
12 & 1264&	197532&	273.95&	34952&	19.53&	293.48&	76&	169814&	259.89&	1&	0.39&	260.28\\
13 & 468&	2602844&	2452.91&	98204&	48.79&	2501.70&	112&	1311943&	1454.68&	1&	0.20&	1454.88\\
14& 299&	3958319&	2835.65&	112485&	98.08&	2933.73&	187&	1755735&	1666.87&	1&	0.30&	1667.16\\ \hline
Average & 443&	810572&	532.79&	28279&	17.27&	\textbf{550.06}&	54&	397298&	321.55&	1&	0.08&	\textbf{321.63}\\
\hline
\multicolumn{7}{l}{add. user cuts: additional user cuts} \\
\hline
\end{tabular} }
\end{table}

\subsection{Application to Power System Networks}\label{section : power system extention}
Reliability concerns have been one of the main focuses of researchers in power system networks since 1960s \citep{billinton1968transmission, nazemi2014modelling}, and has been constantly the topic of recent studies from various aspects \citep{panteli2017metrics, fathabad2020data}. Depending on the location of generators in a power system, the nodes of a power system network can be categorized into generation nodes, load nodes, and hub nodes, i.e., main substations used to transfer power to various points of the system. From the service reliability point of view, the generation nodes are the most important ones since a failure in generation nodes interrupts the supply and results in a blackout at load points. Also, hub nodes are of special interest since they are larger stations to distribute power in smaller areas. The rest of the nodes are mainly load nodes which are not as important as the other groups from the reliability point of view because the power outage usually occurs due to the failure in generation nodes or hub nodes that are used for power distribution. As result, one may be interested in enhancing the power system reliability while considering distributed attacks, e.g., letting the possibility of failure only for the generation and hub nodes. 

\begin{figure}
    \centering
    \includegraphics[width=0.6\textwidth, height=0.34\textheight]{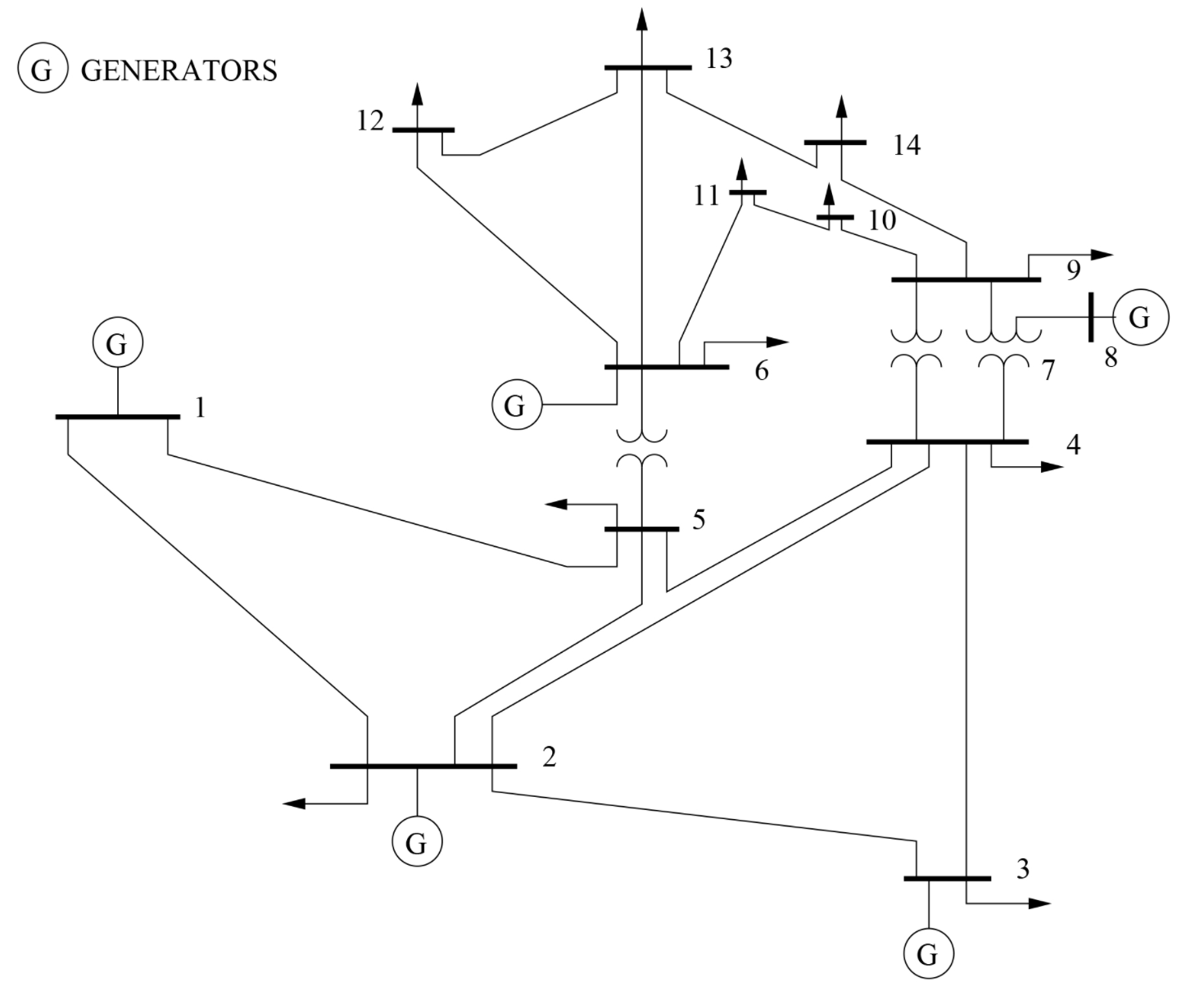}
    \caption{IEEE 14-Bus Test System \citep{christie2000power}}
\label{14-bus}
\end{figure}

\begin{figure}
    \centering
    \includegraphics[width=0.8\textwidth, height=0.33\textheight]{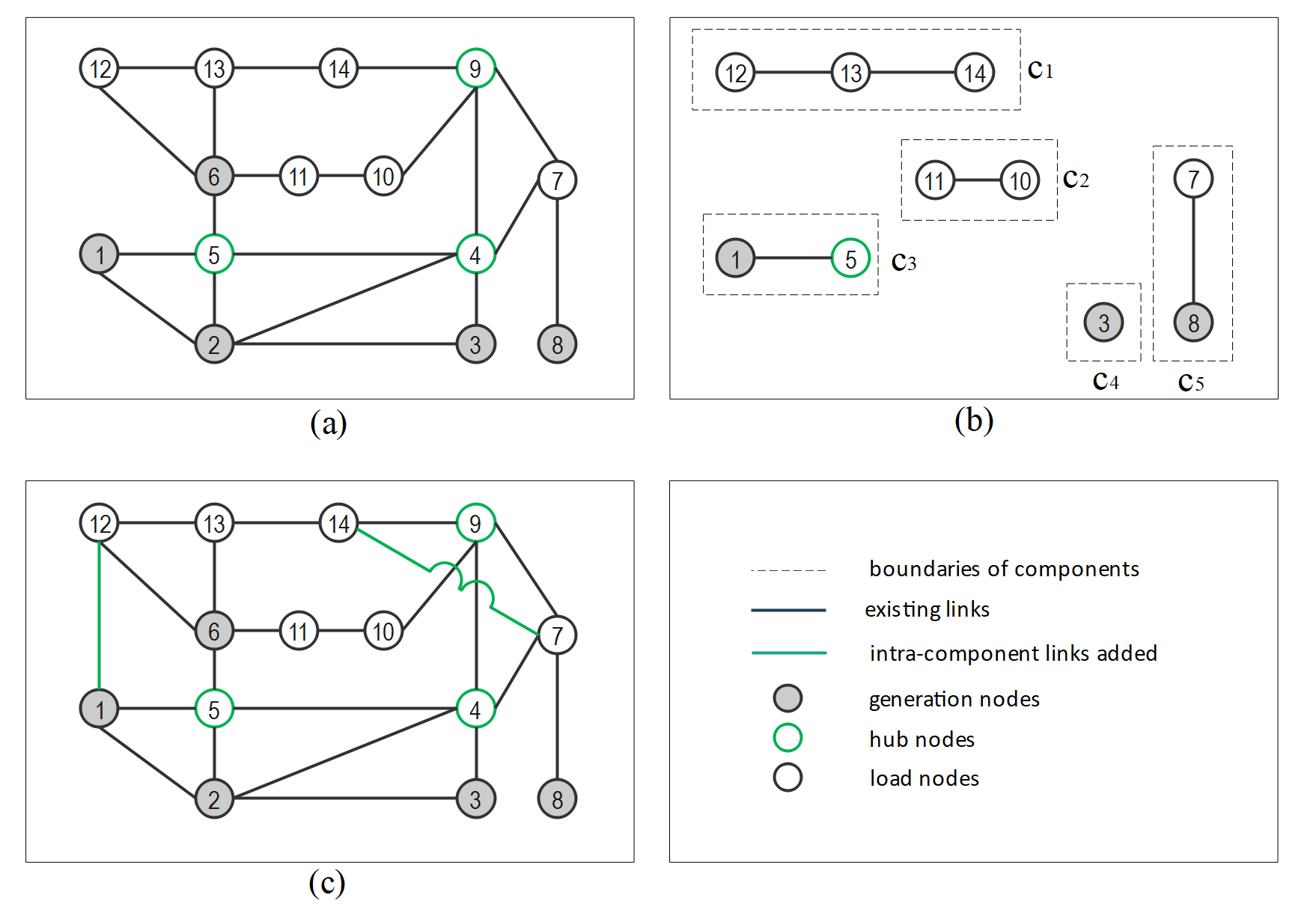}
    \caption{The topology of (a) the underlying network of the IEEE 14 Bus Test system, (b) the attacked network (stage I), (c) the reconstructed network for $B^R=3$ (Stage II)}
 \label{power example}
\end{figure}

In this section, we demonstrate the applicability of the proposed approach in the IEEE 14 Bus Test system representing a portion of the American Electric Power System in the Midwestern US \citep{christie2000power}. Figures \ref{14-bus} and \ref{power example}.a show the actual IEEE 14-Bus system and the underlying graph for the system, where we categorize the buses in 3 groups, namely (i) $\mathcal{N}_g=\{1,2,3,6,8\}$, corresponding to generator nodes, (ii)  $\mathcal{N}_h=\{4,5,9\}$, corresponding to hub nodes, and (iii)  $\mathcal{N}_i=\{7,10,11,12,13,14\}$, corresponding to the rest of nodes that are considered intact in this example.
In such power systems, we consider distributed attacks discussed in section \ref{distributed attacks}, where in the attack problem nodes in $\mathcal{N}_i$ are intact.  As a result, Constraint \ref{s1 c1} is replaced by the following constraints 
\begin{subequations} \begin{align} & \sum_{c \in \mathcal{C}} v^A_{ic} \le 1 & \forall i \in \mathcal{N}_g \cup \mathcal{N}_h \label{dist s1 c1} \\ & \sum_{c \in \mathcal{C}} v^A_{ic} = 1 & \forall i \in \mathcal{N}_i \label{dist s2 c1} \end{align} \end{subequations}
to ensure that all vertices of subset $\mathcal{N}_i$ remain intact since they are not of the attacker's interest. Figure  \ref{power example}.b corresponds to the solution of the attack problem where the worst cut set is $|X|=\{2,4,6.9\}$ corresponding to the minimum resilience of $r^*(G)=1$. It can be observed that components such as  $c_1$ and $c_2$ only contains load nodes while some components may contain generation nodes as well as other types of nodes.

In the response stage, one may consider the technical characteristics of the power system prior to adding links. In particular, one can argue that even though linking components $c_1$ and $c_2$ increases the network resilience due to a higher level of connectivity, this connection might not be beneficial from the service reliability point of view. It is because all nodes of components $c_1$ and $c_2$ are load nodes, and such load nodes cannot be supplied in the absence of $|X|=\{2,4,6.9\}$. Therefore, the connection of components $c_1$ and $c_2$ in the response stage is reasonable only if other components with generation nodes are also connected to components $c_1$ or $c_2$. Let $\mathcal{C}_l$ be the set of components whose nodes are all load nodes and $\mathcal{C}_g$ be the set of components that have at least 1 generator. To impose such technical constraints, we extend the response problem for power system networks by adding the following constraint.
\begin{equation} \label{constraint add power}
    \hat{x}_{\sigma(m,n)} \le \sum_{i \in \mathcal{C}_g} \big(\hat{x}_{\sigma(m,i)} + \hat{x}_{\sigma(n,i)}  \big) \qquad \forall m,n \in \mathcal{C}_l 
\end{equation}
where $\hat{x}_{\sigma(m,n)}$ is a binary variable indicating if a MCEIC link between components $m$ an $n$ is added. Constraint \ref{constraint add power} enforces that if a MCEIC link connects load components $m,n \in \mathcal{C}_l $, there should also exists at least one MCEIC link connecting a component with generators to one of the load components $m$ or $n$. As a result, the power system resilience is optimized with respect to power reliability indices as well rather than merely considering the graph topological properties.

In addition, in the response stage, various indices can be used to calculate matrix $\underline{\mathbf{d}}$ for MCEIC to be considered in the response problem. For example, one may consider the installation cost of MCEIC links as well as the reduced interrupted cost as a resulting of connecting two components. For example from Figure \ref{power example}.b, once load component $c_1$ is connected to component $c_3$ with generator, it is expected that generator located at node 1 transfers power to supply nodes of component $c_1$. Thus, this link addition is more beneficial, and as a result one can consider a low value in matrix $\underline{\mathbf{d}}$ for the MCEIC link connecting components $c_1$ and $c_3$. We form matrix $\underline{\mathbf{d}}$ as the installation cost minus the potential reduction in the interrupted cost of the load nodes. In other words, matrix $\underline{\mathbf{d}}$ represents the financial benefits of the MCEIC links. For example, link $e_1$=$12$-$1$ is more beneficial than $e_2$=$1$-$11$ because $e_1$ can potentially lead to supplying 3 load nodes of $c_1$, i.e., nodes 12, 13, 14, and reduce the interruption cost more. However, $e_2$ can reduce the interruption cost of only two load nodes located in $c_2$, i.e., nodes 10, 11. 
Matrix $\underline{\mathbf{d}}$ used in the second stage is formed as follows:
\begin{equation} \label{dmin power example}
\underline{\mathbf{d}} = \begin{blockarray}{cccccc}
c_1 & c_2 & c_3 & c_4 & c_5 \\
\begin{block}{(ccccc)c}
  - & 2.8 (14\textendash10) & 1.4 (12\textendash1) & 1.7 (14\textendash3) & 1.5 (14\textendash7) & c_1 \\
   & - & 1.8 (1\textendash11) & 1.8 (10\textendash3) & 1.7 (11\textendash8) & c_2 \\
   &  & - & 2.1 (5\textendash3) & 1.8 (5\textendash7) & c_3 \\
   &  &  & - & 1.7 (3\textendash8) & c_4 \\
   &  & &  & - & c_5 \\
\end{block}
\end{blockarray}
\nonumber
\end{equation}
Figure \ref{power example} corresponds to the reconstructed network given the available budget for the link addition is 3 units, i.e., $B^R=3$. It is observed that two links can be added in the reconstructed network. Table \ref{Table: power example} summarizes the results for the IEEE 14 Bus test system. It is observed that $X_1=\{2,4,6,9\}$ correspond to the lowest resilience, i.e., $-r^*_{X_1}(G)=1$, in the initial network before link addition. Once, links $e_1$ and $e_2$ are added to the network, the network resilience under the same cut set $X_1$ enhances to $-r^R_{X_1}(G)=7$, implying the system is more resilient under attacks. Furthermore, in the reconstructed network, set $X_2=\{1,2,6,9\} \subset \mathcal{N}_h \cup \mathcal{N}_g$ becomes the new worst cut set as the dependencies in the reconstructed network change after link addition. Thus, if the attacker has the ability to dynamically respond to our link addition actions, the attacker may choose to disable vertices of $X_2$ resulting in the lowest resilience in the reconstructed network, i.e., $-r^{R*}_{X_2}(G)=5$. However, it is obvious that the reconstructed network is significantly more resilient even under the new (dynamic) worst-case attack, compared to the initial network. The mentioned observations validate the benefits of the proposed solution for enhancing networks resilience and robustness in power system networks.

\begin{table}[!t]
\centering
\footnotesize
\caption{Resilience and Robustness analysis of the IEEE 14 Bust Test System for $B^R=3$}
\label{Table: power example}
\scalebox{0.9}{
\begin{tabular}{lcccc}
\hline
& \multicolumn{2}{c}{Initial Network}&\multicolumn{2}{c}{Reconstructed Network}   \\ \hline
 Cut set & \# of components & $\text{NR}: -r_X(G)$ & \# of components & $\text{NR}: -r^R_X(G)$ \\
\hline
$X_1=\{2,4,6,9\}$ & 5 & \textbf{ 1} $(r^*)$ & 3 & 7 \\ 
$X_2=\{1,2,6,9\}$ & 4 & 2 & 3 & \textbf{5} $(r^{R*})$ \\
\hline
\multicolumn{3}{l}{NR: Network Resilience} \\
\multicolumn{3}{l}{*: optimal objective value of the attack problem} \\
\hline
\end{tabular} }
\end{table}

\section{Conclusion}\label{section: conclusion} 
This paper proposed a novel mixed-integer programming (MIP) formulation for network resilience and robustness under targeted, random, and distributed attacks. The proposed model MIP formulation was derived based on graph-theoretical aspects to increase the network hardiness against attacks, i.e., network robustness, and to maximize the network functionality after attacks, i.e., network resilience. The proposed MIP model consists of two stages to first identify the worst-case attack and then to maximize the network resilience and robustness under the worst-case attack. In addition to the novel MIP model, we showed how the proposed solution method can solve the proposed model efficiently. Particularly, we exploited the structure of the model, derived a tight relaxation, and implemented recent lifted cover inequalities (LCI) in the model. The numerical results showed promising results in both the robustness and computational efficiency of the proposed approach. 

The applicability of the proposed approach was shown in power system networks while considering the technical properties of such networks. A potential research direction is to apply the proposed approach for more applications such as integrated communication and supply chain networks and include their specific problem-related limitations. In addition, it is an interesting future research to derive other cutting plane algorithms for the proposed approach to enhance the computational efficiencies even more.

\bibliographystyle{informs2014}
\bibliography{sample}

\begin{thebibliography}{62}
\providecommand{\natexlab}[1]{#1}
\providecommand{\url}[1]{\texttt{#1}}
\providecommand{\urlprefix}{URL }

\bibitem[{Abbas et~al.(2017)Abbas, Laszka, \protect\BIBand{}
  Koutsoukos}]{abbas2017improving}
Abbas W, Laszka A, Koutsoukos X (2017) Improving network connectivity and
  robustness using trusted nodes with application to resilient consensus.
  \emph{IEEE Transactions on Control of Network Systems} 5(4):2036--2048.

\bibitem[{Ahmadian et~al.(2020)Ahmadian, Lim, Cho, \protect\BIBand{}
  Bora}]{ahmadian2020quantitative}
Ahmadian N, Lim GJ, Cho J, Bora S (2020) A quantitative approach for assessment
  and improvement of network resilience. \emph{Reliability Engineering \&
  System Safety} 200:106977.

\bibitem[{Albert et~al.(2000)Albert, Jeong, \protect\BIBand{}
  Barab{\'a}si}]{albert2000error}
Albert R, Jeong H, Barab{\'a}si AL (2000) Error and attack tolerance of complex
  networks. \emph{nature} 406(6794):378--382.

\bibitem[{Bachmann et~al.(2020)Bachmann, Bustos-Jim{\'e}nez, \protect\BIBand{}
  Bustos}]{bachmann2020survey}
Bachmann I, Bustos-Jim{\'e}nez J, Bustos B (2020) A survey on frameworks used
  for robustness analysis on interdependent networks. \emph{Complexity} 2020.

\bibitem[{Balas(1971)}]{balas1971intersection}
Balas E (1971) Intersection cuts—a new type of cutting planes for integer
  programming. \emph{Operations Research} 19(1):19--39.

\bibitem[{Balas(1975)}]{balas1975facets}
Balas E (1975) Facets of the knapsack polytope. \emph{Mathematical programming}
  8(1):146--164.

\bibitem[{Barefoot et~al.(1987)Barefoot, Entringer, \protect\BIBand{}
  Swart}]{barefoot1987vulnerability}
Barefoot CA, Entringer R, Swart H (1987) Vulnerability in graphs-a comparative
  survey. \emph{J. Combin. Math. Combin. Comput} 1(38):13--22.

\bibitem[{Bertsimas \protect\BIBand{} Sim(2004)}]{bertsimas2004price}
Bertsimas D, Sim M (2004) The price of robustness. \emph{Operations research}
  52(1):35--53.

\bibitem[{Beygelzimer et~al.(2005)Beygelzimer, Grinstein, Linsker,
  \protect\BIBand{} Rish}]{beygelzimer2005improving}
Beygelzimer A, Grinstein G, Linsker R, Rish I (2005) Improving network
  robustness by edge modification. \emph{Physica A: Statistical Mechanics and
  its Applications} 357(3-4):593--612.

\bibitem[{Billinton \protect\BIBand{}
  Bollinger(1968)}]{billinton1968transmission}
Billinton R, Bollinger KE (1968) Transmission system reliability evaluation
  using markov processes. \emph{IEEE Transactions on power apparatus and
  systems} (2):538--547.

\bibitem[{Bixby et~al.(2004)Bixby, Fenelon, Gu, Rothberg, \protect\BIBand{}
  Wunderling}]{bixby2004mixed}
Bixby RE, Fenelon M, Gu Z, Rothberg E, Wunderling R (2004) Mixed-integer
  programming: A progress report. \emph{The sharpest cut: the impact of Manfred
  Padberg and his work}, 309--325 (SIAM).

\bibitem[{Cao et~al.(2013)Cao, Hong, Du, \protect\BIBand{}
  Zhang}]{cao2013improving}
Cao XB, Hong C, Du WB, Zhang J (2013) Improving the network robustness against
  cascading failures by adding links. \emph{Chaos, Solitons \& Fractals}
  57:35--40.

\bibitem[{Chan \protect\BIBand{} Akoglu(2016)}]{chan2016optimizing}
Chan H, Akoglu L (2016) Optimizing network robustness by edge rewiring: a
  general framework. \emph{Data Mining and Knowledge Discovery}
  30(5):1395--1425.

\bibitem[{Christie(2000)}]{christie2000power}
Christie R (2000) Power systems test case archive. \emph{Electrical Engineering
  dept., University of Washington} 108, \urlprefix\url{http://www2.ee.
  washington.edu/research/pstca/}.

\bibitem[{Chv{\'a}tal(1973)}]{chvatal1973tough}
Chv{\'a}tal V (1973) Tough graphs and hamiltonian circuits. \emph{Discrete
  Mathematics} 5(3):215--228.

\bibitem[{Cook et~al.(1993)Cook, Rutherford, Scarf, \protect\BIBand{}
  Shallcross}]{cook1993implementation}
Cook W, Rutherford T, Scarf HE, Shallcross D (1993) An implementation of the
  generalized basis reduction algorithm for integer programming. \emph{ORSA
  Journal on Computing} 5(2):206--212.

\bibitem[{Cozzens et~al.(1995)Cozzens, Moazzami, \protect\BIBand{}
  Stueckle}]{cozzens1995tenacity}
Cozzens M, Moazzami D, Stueckle S (1995) The tenacity of a graph. \emph{Graph
  Theory, Combinatorics, and Algorithms (Yousef Alavi and Allen Schwenk eds.)
  Wiley, New York} 1111--1112.

\bibitem[{Cui et~al.(2017)Cui, Zhu, Xun, \protect\BIBand{}
  Xia}]{cui2017enhance}
Cui P, Zhu P, Xun P, Xia Z (2017) Enhance the robustness of cyber-physical
  systems by adding interdependency. \emph{2017 IEEE International Conference
  on Intelligence and Security Informatics (ISI)}, 203--203 (IEEE).

\bibitem[{Dong et~al.(2015)Dong, Fang, Tian, \protect\BIBand{}
  Zhang}]{dong2015approaches}
Dong Z, Fang Y, Tian M, Zhang R (2015) Approaches to improve the robustness on
  interdependent networks against cascading failures with load-based model.
  \emph{Modern Physics Letters B} 29(32):1550210.

\bibitem[{Dong et~al.(2020)Dong, Tian, Tang, Li, \protect\BIBand{}
  Lai}]{dong2020improving}
Dong Z, Tian M, Tang R, Li X, Lai J (2020) Improving the robustness of spatial
  networks by link addition: more and dispersed links perform better.
  \emph{Nonlinear Dynamics} 100(3):2287--2298.

\bibitem[{Fathabad et~al.(2020)Fathabad, Cheng, Pan, \protect\BIBand{}
  Qiu}]{fathabad2020data}
Fathabad AM, Cheng J, Pan K, Qiu F (2020) Data-driven planning for renewable
  distributed generation integration. \emph{IEEE Transactions on Power Systems}
  35(6):4357--4368.

\bibitem[{Filabadi \protect\BIBand{} Azad(2020)}]{filabadi2020robust}
Filabadi MD, Azad SP (2020) Robust optimisation framework for {SCED} problem in
  mixed {AC-HVDC} power systems with wind uncertainty. \emph{IET Renewable
  Power Generation} 14(14):2563--2572.

\bibitem[{Filabadi \protect\BIBand{}
  Mahmoudzadeh(2021)}]{filabadi2019effective}
Filabadi MD, Mahmoudzadeh H (2021) Effective budget of uncertainty for classes
  of robust optimization. \emph{INFORMS Journal on Optimization} (in press).

\bibitem[{Geoffrion(1974)}]{geoffrion1974lagrangean}
Geoffrion AM (1974) Lagrangean relaxation for integer programming.
  \emph{Approaches to integer programming}, 82--114 (Springer).

\bibitem[{Golan et~al.(2020)Golan, Jernegan, \protect\BIBand{}
  Linkov}]{golan2020trends}
Golan MS, Jernegan LH, Linkov I (2020) Trends and applications of resilience
  analytics in supply chain modeling: systematic literature review in the
  context of the covid-19 pandemic. \emph{Environment Systems and Decisions}
  40:222--243.

\bibitem[{Gomory(1960)}]{gomory1960algorithm}
Gomory R (1960) An algorithm for the mixed integer problem. Technical report,
  RAND CORP SANTA MONICA CA.

\bibitem[{Guan et~al.(2011)Guan, Chen, \protect\BIBand{}
  Qian}]{guan2011routing}
Guan ZH, Chen L, Qian TH (2011) Routing in scale-free networks based on
  expanding betweenness centrality. \emph{Physica A: Statistical mechanics and
  its applications} 390(6):1131--1138.

\bibitem[{Guignard \protect\BIBand{} Spielberg(1981)}]{guignard1981logical}
Guignard M, Spielberg K (1981) Logical reduction methods in zero-one
  programming—minimal preferred variables. \emph{Operations Research}
  29(1):49--74.

\bibitem[{Henry \protect\BIBand{} Ramirez-Marquez(2016)}]{henry2016impacts}
Henry D, Ramirez-Marquez JE (2016) On the impacts of power outages during
  hurricane sandy—a resilience-based analysis. \emph{Systems Engineering}
  19(1):59--75.

\bibitem[{Holme et~al.(2002)Holme, Kim, Yoon, \protect\BIBand{}
  Han}]{holme2002attack}
Holme P, Kim BJ, Yoon CN, Han SK (2002) Attack vulnerability of complex
  networks. \emph{Physical review E} 65(5):056109.

\bibitem[{Ivanov(2020)}]{ivanov2020predicting}
Ivanov D (2020) Predicting the impacts of epidemic outbreaks on global supply
  chains: A simulation-based analysis on the coronavirus outbreak
  (covid-19/sars-cov-2) case. \emph{Transportation Research Part E: Logistics
  and Transportation Review} 136:101922.

\bibitem[{Ji et~al.(2016)Ji, Wang, Liu, Chen, Tang, Wei, \protect\BIBand{}
  Tu}]{ji2016improving}
Ji X, Wang B, Liu D, Chen G, Tang F, Wei D, Tu L (2016) Improving
  interdependent networks robustness by adding connectivity links.
  \emph{Physica A: Statistical Mechanics and its Applications} 444:9--19.

\bibitem[{Jung(1978)}]{jung1978class}
Jung HA (1978) On a class of posets and the corresponding comparability graphs.
  \emph{Journal of Combinatorial Theory, Series B} 24(2):125--133.

\bibitem[{Kannan \protect\BIBand{} Monma(1978)}]{kannan1978computational}
Kannan R, Monma CL (1978) On the computational complexity of integer
  programming problems. \emph{Optimization and Operations Research}, 161--172
  (Springer).

\bibitem[{Kazawa \protect\BIBand{} Tsugawa(2020)}]{kazawa2020effectiveness}
Kazawa Y, Tsugawa S (2020) Effectiveness of link-addition strategies for
  improving the robustness of both multiplex and interdependent networks.
  \emph{Physica A: Statistical Mechanics and its Applications} 545:123586.

\bibitem[{Latora \protect\BIBand{} Marchiori(2007)}]{latora2007measure}
Latora V, Marchiori M (2007) A measure of centrality based on network
  efficiency. \emph{New Journal of Physics} 9(6):188.

\bibitem[{Letchford \protect\BIBand{} Souli(2019)}]{letchford2019lifted}
Letchford AN, Souli G (2019) On lifted cover inequalities: A new lifting
  procedure with unusual properties. \emph{Operations Research Letters}
  47(2):83--87.

\bibitem[{Li \protect\BIBand{} Zhang(2010)}]{li2010extremal}
Li Y, Zhang S (2010) Extremal graphs with given order and the rupture degree.
  \emph{Computers \& Mathematics with Applications} 60(6):1706--1710.

\bibitem[{Li et~al.(2005)Li, Zhang, \protect\BIBand{} Li}]{li2005rupture}
Li Y, Zhang S, Li X (2005) Rupture degree of graphs. \emph{International
  Journal of Computer Mathematics} 82(7):793--803.

\bibitem[{MacKenzie et~al.(2012)MacKenzie, Santos, \protect\BIBand{}
  Barker}]{mackenzie2012measuring}
MacKenzie CA, Santos JR, Barker K (2012) Measuring changes in international
  production from a disruption: Case study of the japanese earthquake and
  tsunami. \emph{International Journal of Production Economics}
  138(2):293--302.

\bibitem[{Matisziw \protect\BIBand{} Murray(2009)}]{matisziw2009modeling}
Matisziw TC, Murray AT (2009) Modeling s--t path availability to support
  disaster vulnerability assessment of network infrastructure. \emph{Computers
  \& Operations Research} 36(1):16--26.

\bibitem[{MCSP. \protect\BIBand{} Payne(2005)}]{bsc2005relaxation}
MCSP BP, Payne RA (2005) \emph{Relaxation techniques} (Elsevier).

\bibitem[{Mens et~al.(2011)Mens, Klijn, de~Bruijn, \protect\BIBand{} van
  Beek}]{mens2011meaning}
Mens MJ, Klijn F, de~Bruijn KM, van Beek E (2011) The meaning of system
  robustness for flood risk management. \emph{Environmental science \& policy}
  14(8):1121--1131.

\bibitem[{Nazemi \protect\BIBand{} Mashayekhi(2014)}]{nazemi2014modelling}
Nazemi A, Mashayekhi M (2014) Modelling welfare loss in iranian electricity
  market. \emph{Energy \& the Economy, 37th IAEE International Conference, June
  15-18, 2014} (International Association for Energy Economics).

\bibitem[{Nemhauser \protect\BIBand{} Wolsey(1990)}]{nemhauser1990recursive}
Nemhauser GL, Wolsey LA (1990) A recursive procedure to generate all cuts for
  0--1 mixed integer programs. \emph{Mathematical Programming} 46(1):379--390.

\bibitem[{Nguyen et~al.(2013)Nguyen, Shen, \protect\BIBand{}
  Thai}]{nguyen2013detecting}
Nguyen DT, Shen Y, Thai MT (2013) Detecting critical nodes in interdependent
  power networks for vulnerability assessment. \emph{IEEE Transactions on Smart
  Grid} 4(1):151--159.

\bibitem[{Panteli et~al.(2017)Panteli, Mancarella, Trakas, Kyriakides,
  \protect\BIBand{} Hatziargyriou}]{panteli2017metrics}
Panteli M, Mancarella P, Trakas DN, Kyriakides E, Hatziargyriou ND (2017)
  Metrics and quantification of operational and infrastructure resilience in
  power systems. \emph{IEEE Transactions on Power Systems} 32(6):4732--4742.

\bibitem[{Papadimitriou(1981)}]{papadimitriou1981complexity}
Papadimitriou CH (1981) On the complexity of integer programming. \emph{Journal
  of the ACM (JACM)} 28(4):765--768.

\bibitem[{Ruj \protect\BIBand{} Pal(2014)}]{ruj2014analyzing}
Ruj S, Pal A (2014) Analyzing cascading failures in smart grids under random
  and targeted attacks. \emph{2014 IEEE 28th International Conference on
  Advanced Information Networking and Applications}, 226--233 (IEEE).

\bibitem[{Sadrabadi et~al.(2021)Sadrabadi, Jafari-Nodoushan, \protect\BIBand{}
  Bozorgi-Amiri}]{sadrabadi2021resilient}
Sadrabadi MHD, Jafari-Nodoushan A, Bozorgi-Amiri A (2021) Resilient supply
  chain under risks: A network and structural perspective. \emph{Iranian
  Journal of Management Studies (IJMS)} 14(4):735--760.

\bibitem[{Scellato et~al.(2011)Scellato, Leontiadis, Mascolo, Basu,
  \protect\BIBand{} Zafer}]{scellato2011evaluating}
Scellato S, Leontiadis I, Mascolo C, Basu P, Zafer M (2011) Evaluating temporal
  robustness of mobile networks. \emph{IEEE Transactions on Mobile Computing}
  12(1):105--117.

\bibitem[{Sen et~al.(2014)Sen, Mazumder, Banerjee, Das, \protect\BIBand{}
  Compton}]{sen2014identification}
Sen A, Mazumder A, Banerjee J, Das A, Compton R (2014) Identification of k most
  vulnerable nodes in multi-layered network using a new model of
  interdependency. \emph{2014 IEEE Conference on Computer Communications
  Workshops (INFOCOM WKSHPS)}, 831--836 (IEEE).

\bibitem[{Shao et~al.(2011)Shao, Buldyrev, Havlin, \protect\BIBand{}
  Stanley}]{shao2011cascade}
Shao J, Buldyrev SV, Havlin S, Stanley HE (2011) Cascade of failures in coupled
  network systems with multiple support-dependence relations. \emph{Physical
  Review E} 83(3):036116.

\bibitem[{Sun \protect\BIBand{} Shayman(2007)}]{sun2007pairwise}
Sun F, Shayman MA (2007) On pairwise connectivity of wireless multihop
  networks. \emph{International Journal of Security and Networks}
  2(1-2):37--49.

\bibitem[{Tits et~al.(2006)Tits, Absil, \protect\BIBand{}
  Woessner}]{tits2006constraint}
Tits AL, Absil PA, Woessner WP (2006) Constraint reduction for linear programs
  with many inequality constraints. \emph{SIAM Journal on Optimization}
  17(1):119--146.

\bibitem[{Wang \protect\BIBand{} Van~Mieghem(2008)}]{wang2008algebraic}
Wang H, Van~Mieghem P (2008) Algebraic connectivity optimization via link
  addition. \emph{Proceedings of the 3rd International Conference on
  Bio-Inspired Models of Network, Information and Computing Sytems}, 1--8.

\bibitem[{Wang(2020)}]{wang2020multiple}
Wang Q (2020) Multiple-attribute decision making-based optimization of
  robustness against cascading failures in interdependent network. \emph{2020
  2nd International Conference on Information Technology and Computer
  Application (ITCA)}, 452--456 (IEEE).

\bibitem[{Wang et~al.(2018)Wang, Zhou, Li, Cao, \protect\BIBand{}
  Lin}]{wang2018improving}
Wang X, Zhou W, Li R, Cao J, Lin X (2018) Improving robustness of
  interdependent networks by a new coupling strategy. \emph{Physica A:
  Statistical Mechanics and its Applications} 492:1075--1080.

\bibitem[{Wolsey(1975)}]{wolsey1975faces}
Wolsey LA (1975) Faces for a linear inequality in 0--1 variables.
  \emph{Mathematical Programming} 8(1):165--178.

\bibitem[{Wolsey \protect\BIBand{} Nemhauser(1999)}]{wolsey1999integer}
Wolsey LA, Nemhauser GL (1999) \emph{Integer and combinatorial optimization},
  volume~55 (John Wiley \& Sons).

\bibitem[{Zhao \protect\BIBand{} Xu(2009)}]{zhao2009enhancing}
Zhao J, Xu K (2009) Enhancing the robustness of scale-free networks.
  \emph{Journal of Physics A: Mathematical and Theoretical} 42(19):195003.

\bibitem[{Zhao et~al.(2015)Zhao, Zhang, \protect\BIBand{}
  Yang}]{zhao2015cascading}
Zhao Z, Zhang P, Yang H (2015) Cascading failures in interconnected networks
  with dynamical redistribution of loads. \emph{Physica A: Statistical
  Mechanics and its Applications} 433:204--210.

\end{thebibliography}

\end{document}